\magnification=\magstep1
                                                                                
\newcount\sec \sec=0
\input Ref.macros
\input math.macros
\input labelfig.tex
\forwardreferencetrue
\citationgenerationtrue
\initialeqmacro
\sectionnumberstrue
\def\vh{{\cal V}(H)}
\def\pk{{\cal P}_k}
\def\pii{\pi}
\def\p_i{{\cal P}_i}
\def\pj{{\cal P}_j}
\def\rk{{\cal R}_k}
\def\ri{{\cal R}_i}
\def\rj{{\cal R}_j}
\def\qi{{\cal Q}_i}
\def\qj{{\cal Q}_j}
\def\Vol{{\rm Vol}}
\def\Area{{\rm Vol_{d-1}}}
\def\dist{{\rm dist}}
\def\max{{\rm max}}

\def\sur{{\rm sur}}

\def\y{y}
\def\b{b}

\def\capa{{\rm cap}}
\input epsf

\title{Invariant matchings of exponential tail on coin flips in $\Z^d$}
\author{\'Ad\'am Tim\'ar}
\bigskip

\abstract{Consider Bernoulli(1/2) percolation on $\Z^d$, and define a 
perfect matching between open and closed vertices in a way that is a 
deterministic 
equivariant function of the configuration. We want to find such matching 
rules that make the probability that the pair of the origin is at distance 
greater than $r$ decay as fast as possible. For two dimensions, we give a 
matching of decay $cr^{1/2}$, which is optimal. 
For dimension at least 3 we give a matching rule that has an exponential 
tail. This substantially improves previous bounds. The construction has 
two major parts: first we define a sequence of 
coarser and coarser partitions of $\Z^d$ in an equivariant way, such 
that with high 
probability the cell of a fixed point is like a cube, and the labels in it 
are i.i.d. Then we define a matching for a fixed finite cell, which 
stabilizes as we repeatedly apply it for the cells of the consecutive 
partitions. Our methods also work in the case when one wants to match 
points 
of two Poisson processes, and they may be applied to allocation 
questions.  
}
                                                                                
\bottomII{Primary 60K35, 82B43. Secondary 60B99.}
{Point processes, invariant matching, optimal matching.}
{Research partially supported by Hungarian National Foundation for 
Scientific Research Grant TO34475.}
                                    
\bsection{Introduction}{s.intro}

Fix $(\Omega, \Sigma, \P)$, where $\Omega=\{0,1\}^{Z^d}$, $\Sigma$ is the 
product $\sigma$-algebra, and $\P$ is the product of Bernoulli 
measures with parameter $1/2$. 
We prove the following theorems.

\procl t.dim2
For $d= 1,2$, there exists a deterministic perfect matching 
$\phi_\omega=\phi$ between 
$\{ 
x\in 
\Z^d\, :\, 
\omega (x)=0\}$ and $\{ x\in \Z^d\, :\,
\omega (x)=1\}$, such that for almost every $\omega\in\Omega$, $\phi$ is an equivariant function (i.e., 
$\phi_\omega =\phi_{g(\omega)}$ for every translation $g$ of $\Z^d$), and
for any $r>0$,
$$\P[\dist (o,\phi (o)) >r]<{c\over 
r^{d/2}}$$
with some constant $c$. 
\endprocl

\procl t.dim3
Consider $d\geq 3$, and $\epsilon>0$ arbitrary. Then there exists a 
deterministic perfect matching
$\phi_\omega=\phi$ between
$\{
x\in
\Z^d\, :\,
\omega (x)=0\}$ and $\{ x\in \Z^d\, :\,
\omega (x)=1\}$, such that for almost every $\omega\in\Omega$, $\phi$ is an equivariant function, and
for any $r>0$,
$$\P[\dist (o,\phi (o)) >r]< C\exp(-c r^{d-2-\epsilon}) .$$
\endprocl

The bound in \ref t.dim3/ can be slightly tightened, see \ref r.f/.
We also have some ideas that could possibly remove the ``$-\epsilon$" from the 
bound, but at the cost of much extra complication. However,
the correct magnitude of the exponent is not known: the only lower bound 
is the trivial $C\exp(-cr^d)$.

For $d=2$ \ref t.dim2/ is new; the best known result has been that of \ref
b.S/.
Note that for $d=1,2$, \ref t.dim2/ is essentially tight by a
theorem of \ref b.HP/, which says that for any matching rule $\phi$, $\E[\dist
(o,\phi(o))^{d/2}]=\infty$ for these dimensions.
For higher dimensions it was believed by Holroyd and Peres that there
would be an exponential
bound (see \ref b.HP/ and also \ref b.S/).

Our proofs rely on the following theorem, which is of independent 
interest. Informally, it claims that there is a sequence of coarser and 
coarser partitions for the space that 
are deterministic functions of the point 
configuration, and still most of the cells are (approximate) cubes with 
i.i.d. Bernoulli labels in them. 

A subset $X$ of $\Z^d$ ($\R^d$) is called a $k$-{\bf pseudocube}, if $X$
contains some $[k/2]\times\ldots\times [k/2]$ cube, and $X$ is contained
in
some
$2k\times\ldots\times 2k$ cube (these two are referred to as the {\bf
volume condition}), and finally, if $|\partial X|\leq c_0
|X|^{{d-1\over
d}}$, with $c_0=2^{2d}+2^d d$ (referred to as the {\bf isoperimetry 
condition}). This choice for $c_0$ is rather arbitrary, and any greater constant could be used; in particular, 
it is clear that the intersection of $\Z^d$ with a pseudocube of $\R^d$ is a 
pseudocube itself, with a $c_0$ that is only worse by some constant factor.

\procl t.main
Fix $o\in\Z^d$.
There exists a sequence $\{ {\cal Q}_i\}$ of coarser and coarser
partitions of $\Z^d$ that are
equivariant functions of the configuration $\omega\in \Omega$, and such
that there is an event $A_i$ with
$\P [A_i]\geq 1-c2^{-2^{i-1}}$, such that conditioned on
$A_i$, the following hold:

\item{(i)} For each $j\geq i$, the cell $C_j$ of $o$ in ${\cal Q}_j$ is a
$2^j$-pseudocube;
\item{(ii)} conditioned further on the location of $C_j$, the labels in 
$C_j$ are
i.i.d. Bernoulli($1/2$).
\item{(iii)} There is some infinite subsequence ${\cal Q}_{\alpha (k)}$
such
that if $i=\alpha (k)$ then
$C_i$ is a $2^i$ by $2^i$ cube.
\endprocl

In fact the probability of $A_i$ with the above properties can be made
arbitrarily large by an appropriate choice of the parameters in our
construction.

\ref t.dim2/ is a straightforward corollary of \ref t.main/.

\proofof t.dim2
Given the sequence of partitions, the matching is defined similarly to
\ref b.S/.
Namely, consider the sequence ${\cal Q}_i$ from \ref t.main/, and as
$i=1,2,\ldots$, for each $2^i$-pseudocube, match as many yet unmatched 
points as
possible, each with a point of opposite label, but otherwise arbitrarily. 
The central limit theorem
gives
the claim.\Qed

A question similar to the ones above is when one considers a Poisson point 
process (the natural 
generalization of a set of uniformly distributed points in the unit cube 
to an infinite domain), and colors each of the configuration points 
independently red 
or blue with 
probability 1/2. This is the same as taking two independent Poisson point 
processes of the same intensity. Our goal in this setting again is to give 
an ``optimal" perfect matching between the red and blue points, by some 
matching rule that is a deterministic and equivariant function of the 
random point set. (Informally, the matching rule is defined using the 
locations and colors of the configuration points, but no background 
information from the underlying space.) By an optimal matching we meant 
in the $\Z^d$ case that the function $F(r):=\P[$distance of 0 from 
its pair is
$>r]$ decays as fast as possible. For the Poisson case, similarly, we 
condition on that 0 is a configuration point, and want to make 
the tail $F(r):=\P[$distance of 0 from its pair is
$>r\, |\, 0\in\omega]$ tend to 0 fast. The setting of our question shows 
that requiring 
the matching rule to be invariant is natural, since it essentially means 
that $F(r)$ does not change if we replace 0 by any other point. 

We will phrase and prove our 
theorems for the 
$\Z^d$ case, but
our methods can be easily adjusted to the matching problem for 
Poisson point processes. 
Furthermore, a sequence of partitions as in \ref t.main/ can be obtained 
as an equivariant function of {\it one} Poisson point process, by the 
natural modifications if our proof. Because of possible applications, 
let us state this separately. 

\procl t.main2
Fix $o\in\R^d$.
There exists a sequence $\{ {\cal Q}_i\}$ of coarser and coarser
partitions of $\R^d$ that are
equivariant functions of the configuration $\omega$ of a Poisson point 
process on $\R^d$, and such
that there is an event $A_i$ with
$\P [A_i]\geq 1-c2^{-2^{i-1}}$, such that conditioned on
$A_i$, the following hold:

\item{(i)} For each $j\geq i$, the cell $C_j$ of $o$ in ${\cal Q}_j$ is 
a
polyhedron and a $2^j$-pseudocube;
\item{(ii)} for any Lebesgue measurable $A\subset \R^d$, $\P 
\bigl[|\omega\cap A\cap 
C_i|\,\bigl | \,C_i\bigr]={\rm Poisson} (|A\cap C_i|)$.
\item{(iii)} There is some infinite subsequence ${\cal Q}_{\alpha (k)}$
such
that if $i=\alpha (k)$ then
$C_i$ is a $2^i$ by $2^i$ cube.
\endprocl

Here ${\rm Poisson}(\lambda)$ denotes the Poisson distribution of 
intensity 
$\lambda$.

The question that we address was first asked by Holroyd and Peres in \ref 
b.HP/. The best 
results were of order $r^{-s(d)}$, where $s(2)<2/3$, and 
$s(d)<2$ for any $d$, 
\ref b.S/. Similar, but 
considerably simpler matching problems were fully solved in \ref b.HPPS/; 
we will go into the details later in this introduction. Our question is 
related to two intensively studied families of problems. First, invariant 
measurable functions ({\bf factors}) of point processes have been of 
interest from a 
statistical point of view (e.g. Palm processes, allocations; see \ref 
b.Th/, 
\ref b.HL/ and further references therein), from the optimization aspect
(such as minimal spanning trees, \ref b.Ale/),
and from a more general interest about how much information can be 
extracted from a point process (\ref b.HPregi/ being a seminal paper in 
the area). Close relatives to the matching questions treated here are the 
so called {\bf allocation questions}, when one has to 
assign 
disjoint sets of measure 1 to every point of a Poisson point process of 
intensity 1 so that this partition of $\R^d$ is a factor. We will return 
to this later.
The other related branch of problems, from a different direction, is 
finding an optimal matching between independently distributed points in a 
unit 
cube. While this latter field is almost fully explored (as a result of 
work by Ajtai-Koml\'os-Tusn\'ady, Talagrand, Shor, Yukich; see \ref 
b.Y/, \ref b.Ta/ for surveys and a still standing challanging problem by 
Talagrand), the methods there do not seem to apply for our setting, 
because of 
the difficulties arising from the infinite setting and invariance.

While the matching rule we are looking for has to be a {\it deterministic} 
function of the configuration, one may relax this requirement and allow 
extra randomness. The additional freedom we gain this way is considerable. 
An example is that using extra randomness allows one to partition $\R^d$ 
($\Z^d$) to ``nice" subsets (e.g. cubes) with independent labels in it: 
simply take $k\Z^d +v$ with $v\in [0,k]^d$ chosen uniformly. (This 
partition cannot be defined as a factor, by ergodicity.) A partition with 
these properties enables one to use local matching rules, repeated for 
countably many, coarser and coarser partitions. 

A variant setting to our problem is when, instead of matching 
points of 
two colors, we have one Poisson point process, and want to find a perfect 
matching of optimal tail on its points. (For the $\Z^d$ case the analogous 
problem 
is meaningless.) Call this a {\bf 1-matching problem} to distinguish from 
the {\bf 2-matching problem} defined earlier.
Similarly to the matching question in a finite box, \ref 
b.Y/, the 1-matching problem is much simpler than the 2-color case. The 
main reason for 
this is that much of the difficulties in the 2-color problem is 
coming from the difference between the number of vertices of the two 
colors within some given box. This discrepancy is around the square root 
of 
the number of points in the box, and it gives a lower bound to the number 
of points that cannot potentially be matched within the box. 
Since these 
points have to find their pairs beyond the boundary of the box, 
isoperimetry starts playing a role, and this is responsible for the 
dramatic 
change between dimensions $\leq 2$ and $\geq 3$, as seen in \ref t.dim2/ 
and \ref t.dim3/. Obviously, the difficulty coming from discrepancies does 
not arise in the 
1-color case.

Let us summarize briefly, what has been known about the four problems 
given by 1 and 2-color matchings, with randomized or deterministic 
matching schemes. 
See \ref b.HPPS/ for a detailed account and an instructive table.
The randomized 1-color matching has a sharp tail of order $\exp 
(-cr^d)$ for all dimensions. Similarly, for the randomized 2-color problem 
\ref b.HPPS/ obtained a sharp $\leq cr^{-d/2}$ decay for $d=1,2$, and 
$\leq \exp (-cr^d)$ decay for $d\geq 3$. For the 1-color factor matching, 
the tight bounds are of order $\leq c r^{-1}$ for $d=1$, and $\leq \exp 
(-cr^d)$ for $d\geq 2$. For the 2-color factor matching, the best known 
upper bound was $r^{-s(d)}$ with $s(d)\leq 2/(1+4/d)$, \ref b.S/. 

Let us point at another interesting phenomenon about the decay rates of 
various problems. First, there is a big gap between the optimal decay 
rates of the randomized and the deterministic 1-color matching problems 
in dimension 1. 
Since the distance between two neighboring 
configuration points of a Poisson point process in 1 dimension follows an 
exponential distribution, it is rather the slow, linear decay for the 
factor matching that is 
surprising. It sheds some light on how restrictive the requirement of 
giving a {\it deterministic} perfect matching is.
  
Although \ref t.dim3/ is a big improvement to earlier results, the degree 
of the 
optimal rate of decay is still open. A trivial lower bound is the 
following:

\procl l.lower
For any 2-matching scheme on $\Z^d$ 
$$\E[\exp (cr^d)]=\infty,  \label e.pici  
$$
where r is the distance between the origin and its pair, and $c$ is some 
positive constant.\endprocl

Similar statement holds for the Poisson case, as mentioned before. For a 
proof observe that the distance of the origin to the closest point of 
opposite color is a lower bound for $r$, and this already 
satisfies \ref e.pici/.

Though allocation questions have a flavor similar to matching questions, 
we do not know of any direct connection that would make them equivalent in 
some sense. If $\omega_1$ and $\omega_2$ are the configurations of two 
independent Poisson point processes of intensity 1 (call them yellow and 
blue 
points), and $A_i$ is a deterministic, invariant allocation rule 
($i=1,2$), so that $A_i(x)$ is a 
set of measure 1 for each $x\in \omega_i$, and these partition $\R^d$ for 
each $i$, then define a bipartite graph $\Gamma=(V,E)$; 
$V=\omega_1\cup\omega_2$. Namely, let $x_1$ and $x_2$ be adjacent if 
$x_1\in \omega_1$, $x_2\in \omega_2$, and $A(x_1)\cap A(x_2)\not 
=\emptyset$. Then K\"onig's theorem (generalized to locally finite 
infinite graphs) implies that there is a perfect 
matching in $\Gamma$. Hence, as observed by Holroyd and Peres, if one 
could 
define a perfect matching for $\Gamma$ in an {\it invariant} way, then one 
would get a matching rule between $\omega_1$ and $\omega_2$ that has 
essentially the same tail behavior as the allocation rule (at least for 
tails that decay relatively fast). Hence, the existence of an invariant 
perfect matching could be used to give a perfect matching 
from an allocation rule. The best known allocation rule so far is the 
so-called
gravitation allocation, \ref 
b.gravi/.
Although the matching scheme presented in the present paper has a better 
tail 
than
what is proved for the gravitation allocation rule in \ref b.gravi/, 
the following question is still of interest.

\procl g.kerdes
Let $G=(V,E)$ be a random locally finite graph with $V\subset \R^d$, and 
such that for any isometry (translation) $g$ of $\R^d$, $(gV,gE)$ has the 
same distribution as $(V,E)$. Suppose that there is a perfect matching in 
$G$ almost always. Then there is also a perfect matching that is a 
deterministic equivariant function of $G$.\endprocl

The conjecture has a similar flavor to one asked by Bowen and Lyons, 
whether every quasi-transitive planar graph has a periodic 4-coloring, or 
a question of Lyons and Schramm, whether every infinite 
quasi-transitive graph 
has an invariant random coloring with as many colors as its cromatic 
number.

While \ref g.kerdes/ leaves the question of creating  matching schemes 
from allocation 
schemes open, the current proof for the matching question provided us 
with 
a tool for the allocation question. 
In joint work with Ander Holroyd,
we are planning to apply the sequence of partitions in \ref t.main2/ to 
create an allocation rule using the technique in \ref b.AKT/, 
which we
believe may have a tight $\exp (-cr^d)$ tail.

A standard tool for the study of invariant processes is called Mass 
Transport Principle (MTP); see \ref b.HPPS/ for a version that is 
close 
to our setting (Lemma 8), and also for further references.
We will use the MTP via two of its straightforward corollaries, which we 
state separately:

\procl l.MTP 
Suppose there is a given $\epsilon >0$, a ${\cal P}$ invariant partition 
of $\R^d$ ($\Z^d$) to 
measurable 
sets, and an invariant measurable subset $S$ of $\R^d$ ($\Z^d$). If every 
$C\in {\cal P}$ satisfies $|C\cap S|/|C|\leq\epsilon$, then for any point 
$o$ of 
$\R^d$ ($\Z^d$) one has $\P[o\in S]\leq \epsilon$. \endprocl

\procl l.parja
Let $S$ be some random measurable subset of $\R^d$ ($\Z^d$), and $M$ be an 
invariant 
perfect matching on the points of $\omega$. Then 
$\P[o$ or $M(o)$ is in $S]\leq 2\P[o$ is in $S]$. \endprocl

In \ref s.partitions/ we present a sequence of partitions as in \ref 
t.main/. 
In \ref s.matching/ we prove \ref t.dim3/. Namely, we present a 
matching algorithm for Bernoulli labelled points in some {\it fixed cube} 
and then this is used to 
give the desired invariant perfect matching on $\Z^d$. We will 
repeatedly apply the algorithm for bigger and bigger cubes coming from 
the ${\cal Q}_{\alpha (i)}$, using the matching algorithm in a way that we 
only rematch 
vertices of smaller and smaller density. 

\ref s.partitions/ and \ref s.matching/ are independent, except for that 
\ref s.matching/ uses \ref t.main/. Similar ternms and notation may have different definitions in the two sections.

In the rest of the paper $c$ and $c'$ always denote positive constants 
depending 
only on $d$, and 
their values may change from line to line.

\bsection{The sequence of partitions}{s.partitions}

In this section we construct the sequence of partitions of \ref t.main/.
Since we achieve this using consequtive sequences of partitions, a look at 
the summary of \ref r.summary/ may fascilitate the reader. 

We shall think of $\Z^d$ as embedded in $\R^d$. We shall define 
Voronoi tessallations of $\R^d$, and then other partitions based on that, always using some 
subset of
$\Z^d$, chosen as a deterministic equivariant function of the labelling. 
In all these 
cases there is also an inherited partition for $\Z^d$, defined by the 
cells of the tessallation. The reason we prefer to partition $\R^d$ is because it sheds light on how the proof works for Poisson point processes, and also because some geometric arguments are simpler to phrase this way. On the other hand it is clear that if we construct the desired sequence of partitions for $\R^d$, that gives rise to a partition for $\Z^d$ as in \ref t.main/.
For any subset $A$ of $\R^d$, we say that $A$ is 
measurable, if it is Lebesgue measurable. Denote by $|A|$ its Lebesgue 
measure. On $\Z^d$ the $\sigma$-algebra that we consider is the discrete one, 
and $|A|$ stands, as usual, for the number of elements in $A$. 
If ${\cal P}$ is some partition, then the sets that it consists of are called the {\it classes} or {\bf cells} of ${\cal P}$. When $o$ is a point of $\Z^d$, 
the cell that contains $o$ is denoted by ${\cal P}(o)$.

We say that a partition $P$ is a {\bf refinement} of partition $Q$, if any 
two 
elements
in the same class of $P$ are also in the same class of $Q$. If $P$ is a 
refinement of $Q$, then $Q$ is a {\bf coarsening} of $P$.
By the union of partitions $P_i$ we mean their coarsest common refinement, 
that is, 
the partition where two 
elements are in the same class if and only if they are in the same class 
for each $P_i$. We denote this partition by $\vee P_i$.
Finally, if ${\cal P}$ is a partition, ${\cal Q}$ is some set of pairwise disjoint subsets of $\Z^d$, then the finest common coarsening of ${\cal P}$ and ${\cal Q}$ is the partition defined by the equivalence relation where two elements are equivalent if they belong to the same class of ${\cal P}$, or the same set in ${\cal Q}$.

The proof of the next lemma is straightforward.

\procl l.coarse
If ${\cal P}_i'$ is a sequence of coarser and coarser partitions of $\Z^d$, ${\cal Q}$ is a set of pairwise disjoint subsets of $\Z^d$, and ${\cal P}_i$ is the finest common coarsening of ${\cal Q}$ and ${\cal P}_i'$, then ${\cal P}_i$ is a sequence of coarser and coarser partitions.\endprocl


In this section {\bf balls} are understood in 
the infinity norm 
(except for two places, where we refer to the ``usual" ball as the {\it norm-2 
ball}). That is, 
the ball $B(x,r)$ of radius $r$ around a point $x=(x_1,\ldots, x_d)$ is 
$\{(y_1,\ldots,y_d)\in \Z^d\, :\, \max_i |x_i -y_i|\leq r\}$. Similarly, 
the {\bf $r$-neighborhood} of a set $H\subset R^d$ is the union of all 
the balls 
of radius $r$ around some point of $H$. Hence the terms ``cube" and 
``ball" stand for the same objects, unless otherwise mentioned.

We denote by $\Vol (A)$ the volume of a $d$ dimensional polyhedron $A$, 
and by 
$\Vol_{d-1} (\partial A)$ the surface area of $A$. By a {\bf path} in a 
graph we 
always understand a simple path, that is, a (finite) sequence of 
vertices with no 
repetitions, such that any two consecutive ones are connected by an edge.

If $H$ is some discrete subset of $\Z^d$, let $\vh$ be the Voronoi 
tessallation of $\R^d$ determined by $H$. Given a cell of some Voronoi 
tessallation $\vh$, we call the (unique) element of $H$ in the cell the 
{\bf centre} of the cell.

Denote by $S_k$ the set of points $x\in\Z^d$ with the property that any 
vertex $y$ such 
that $||x-y||_\infty=i$, $i\in\{0,1,\ldots,k\}$, satisfies $\omega (y)= 
(-1)^i$. Given $x\in \Z^d$, call the configuration on the set 
$\{y\,:\, ||x-y||_\infty \leq k\}$ a {\bf $k$-bulb of} $x$, 
if $\omega (y)=(-1)^i$ whenever $||x-y||_\infty =i$, $i\in\{0,1,\ldots, 
k\}$, i.e. if $x\in S_k$. 
Thus $S_k$ is defined as an equivariant function of the random 
labelling.
A simple but important consequence of the definition is that 
$\cup_{n=k}^\infty S_n =S_k$ is ``sparse": the 
probability of being in $S_k$ is $2^{-ck^d}$.
Another reason for the choice of $S_k$ is
that any two elements of $S_{k}$ have distance at least $2k$ from 
each other, because the $k$-bulbs of two elements in $S_k$ can intersect 
only in their  
boundaries. Hence in ${\cal V} (S_{k})$ every Voronoi cell contains a 
norm-2 ball 
of radius $k$. The most important property of $k$-bulbs is stated in the 
next lemma.


\procl l.giant 
Let $t\in\Z^+$, $o$ be a vertex of $\Z^d$. Let $B:=B(o,t)$, and for a 
vertex $x\in\Z^d$ let $r_{t,x}:=\max (2t, \dist (x,o))$. Call $x$ 
a $t$-giant (of the configuration $\omega$), if the 
configuration on $B(x,r_{t,x})\setminus B$ inherited from $\omega$ can 
be extended to $B(x,r_{t,x})$ so that we get an $r_{t,x}$-bulb 
around $x$. Then
$$\P [\exists\; t{\rm -giant}]\leq C c^{-t^d}$$
with some $C,c>1$.\endprocl

\proof 
Clearly $|B(x,r_{t,x})\setminus B|\geq {|B(x,r_{t,x})|\over 2}$. Hence the 
following hold:
$$\E[{\rm number \; of \;} t{\rm -giants}]\leq \sum_{x\in\Z^d} 
2^{-|B(x,r_{t,x})\setminus B|}\leq$$
$$\leq |B| 2^{-|B(x,2t)|/2}+
\sum_{i=t+1}^\infty c_1 
i^{d-1} 2^{-{|(B(x,i)|\over 2}}\leq (2t)^d 
2^{-2^{d-1}t^d}+\sum_{i=t+1}^\infty c_1 i^{d-1} 
c_2^{-i^d}\leq Cc^{-t^d},$$
with some constants $C$ and $c_1,c_2,c>1$. By Markov's inequality the
same upper bound is valid for the probability that there exists a 
$t$-giant. \Qed

Given some grid $v+r\Z^d$,$v\in \Z^d,\; r\in \Z$, define 
the {\bf basic cubes of} 
$v+r\Z^d$ to be the cubes that have the form $v+\{(x_1,\ldots,x_d)\in\R^d 
\,:\, r 
a_i\leq x_i < r(a_i+1)\}$ with some integers $a_i$.

For a $C\subset \Z^d$, denote by $\partial_\rho C$ 
the $\rho$-neighborhood of the boundary $\partial C$ of $C$. Fix a point $o\in 
Z^d$. 

\procl l.isoperi Let $K$ be a convex polyhedron in $\R^d$, and $K'$ be the 
union of basic cubes of $r\Z^d$ fully contained in $K$.
If $K$ contains a norm-2 ball of radius $R$ around a point $o$, then 
$${\Vol (K\setminus K')\over \Vol (K)}\leq {cr\over R},$$
with some constant $c$ depending only on the dimension.
\endprocl

\proof
Any point $z$ of $K\setminus K'$ has some point of $\partial K$ in its 
$r$-neighborhood, otherwise any $r\times\ldots \times r$ cube containing 
$z$ is contained in $K$, in particular $z\in K'$. Hence $K\setminus K'$ 
is contained in $\partial_r K$, and
$$\Vol (K\setminus K')\leq \Vol (K\cap \partial_r K)\leq r \Area (\partial 
K).$$ 

The last inequality is true 
by the convexity of $K$.
On the other hand, ${1\over d!} \Area (\partial K)R\leq \Vol (K)$, which 
can be seen by subdividing $K$ to pyramids with apice in $o$ 
(similarly to
\ref b.T/). Putting the two inequalities together gives the statement.\Qed

Given a subcube $H=\sqcap_{i=1}^d [m_i, m_i +t]$ of $\R^d$, call
hyperfaces of the form $\{x\in K\, :\, x_i=m_{i}+t\}$, with some $i$, {\bf
right
faces}. Let ${\cal P}$ be a partition of $\R^d$.
Define $\theta (H)$ {\it with respect to} ${\cal P}$ (${\cal P}$ treated as a 
hidden parameter of $\theta$ for simplicity) to be the union of ${\cal 
P}$-cells 
that are contained in $H$ or intersect only right-faces of $H$. (See 
Figure 1 for an example.)

\procl l.lepcsos
Let $H$ 
be a $2^i\times\ldots\times 2^i$ cube of the above form (i.e., with each 
1-face parallel to some
coordinate axis), $i'<i$, and ${\cal P}$ be a partition of $\R^d$ to 
convex 
cells 
$C$ that satisfy one of the following:
\item{(i)} $C$ is a $2^{i'}\times\ldots\times 2^{i'}$ cube with 
each 1-face 
parallel to some coordinate axis,
\item{(ii)} $C$ has diameter at most $1$.

Then $\theta (H)$ 
satisfies 
$|\partial \theta (H)| \leq c|\partial H|=c'2^{i(d-1)}$ with constants $c,c'$ 
depending only 
on the dimension.
\endprocl

\proof
Note
that $\partial \theta (H)$ is contained in the boundary of some ${\cal
P}$-cells each of which
intersects $\partial H$. These cells are either $2^{i'}\times\ldots\times 
2^{i'}$ cubes (call their
set $S_1$), or cells of diameter $\leq 1$ (call their set $S_2$).   
Let $X\in S_1$ and $\delta(X)=\delta$ be the minimal number such that $X$ 
intersects some $\delta$ dimensional face of $H$. Then $|X\cap H|\geq 2^{i'\delta}$, because the edges of $X$ are parallel to the
coordinate axes.
Since the elements of $S_1$ are disjoint, we get
that
$|S_1|\leq\sum_{\delta=0}^{d-1} {2^{i\delta} \over 2^{{i'}\delta}}
\times$number 
of $\delta$ dimensional hyperfaces of $H$. That is,  $|S_1|\leq c 
2^{(i-i')(d-1)}$ 
with
some $c$ depending only on $d$. We obtain:
$$|\partial \theta (H)|\leq 2d2^{i'(d-1)}|S_1|+|\cup_{X\in S_2} \partial 
X|\leq
c' 2^{i'(d-1)}2^{(i-i')(d-1)}+|\cup_{X\in S_2} \partial X|.$$
For the second term here, we can use the crude upper bound $|\partial
X|\leq
2^d$ for $X\in S_2$, and $|S_2|\leq |\partial H|$. Hence we have
$$|\partial \theta(H)|\leq c2^{i(d-1)}
+2^d |\partial H|\leq c |\partial
H|=c' 2^{i(d-1)}, \label e.cukisag
$$
using that $H$ is a $2^i\times\ldots\times 2^i$ cube.
\Qed

Fix sequences $\{a_i\}$ and $\{b_i\}$ to be
$b_i:=2^{2^i}$ and $a_i:=2^{b_i}$.

Consider the refinement $\pk$ of ${\cal V} (S_{a_k})$ where we 
partition 
each cell $C$ of ${\cal V} (S_{a_k})$ with center $v(C)$ using the basic 
cubes of 
$v(C)+b_k \Z^d$. That is, two points of $\R^d$ are 
in the same class (cell) of $\pk$, if they are in the same cell $C$ of 
${\cal V}(S_{a_k})$ and the same basic cube of $v(C)+b_k \Z^d$.
By 
\ref 
l.isoperi/ and \ref l.MTP/, we obtain:

\procl l.k-bad 
The probability that $o$ is not in a $b_k\times\ldots\times b_k$ cell of 
${\cal P}_k$ is at most $\P[o$ is in 
$\cup_{K\in \pk}\partial_{b_k} K]\leq
c b_k/a_k$. \endprocl

Next, let $\rj '= 
\vee_{i=j}^\infty {\cal P}_{i}$. We will prove that:

\procl l.seged
The probability that $o$ is not in a $b_j\times\ldots\times b_j$ cube of
$\rj '$ is $\leq \sum_{k=j}^\infty cb_k/a_k + \sum_{i=j+1}^\infty c_0 
b_j/b_i \leq c2^{-2^j}$. \endprocl

There are two possible reasons for $\rj '(o)$ not to be a 
$b_j\times\ldots\times b_j$ cube: either the $\pj (o)$ already 
fails to be a $b_j\times\ldots\times b_j$ cube (in this case we say that 
$o$ is $b_j${\bf -bad}),
or $\pj (o)$ is 
intersected by some cell boundary of $\p_i$, $i>j$. The bounds for these 
two are provided by \ref l.k-bad/ and \ref l.kicsi/, and hence \ref 
l.seged/ will follow.

\procl l.kicsi Suppose $k\in \Z^+$, $\rho >0$. The probability that the
ball of radius $\rho$ around $o$
is intersected by some cell-boundary from $\pk$ is at most
$c_0\rho/b_k$, where $c_0$ is some constant
independent of $\rho$ and $k$.\endprocl

\proof
It is clear that for any $b_k\times\ldots\times b_k$ cell $C$ in $\pk$,
$\Vol (\partial_\rho
C)/\Vol (C)\leq c' \rho b_k^{d-1}/b_k^d=c'\rho/b_k$ with some constant 
$c'$,
and thus \ref l.MTP/ shows that the probability that
$B(o,\rho)$ is
intersected by some cell
boundary is at most $\P[o$ is $b_k$-bad$]+c'\rho/b_k\leq c/2^{b_k}
+c'\rho/b_k\leq c_0\rho/b_k$ with some constant $c_0$, also using \ref
l.k-bad/.
\Qed

\proofof l.seged
By \ref l.kicsi/ and \ref l.k-bad/, the probability that $o$ is in a cell
of $\ri '$ that does
not
coincide with a $b_i\times\ldots\times b_i$ cell of ${\cal 
P}_{i}$ is at most $\P[o$ is
$b_i$-bad$]+\P[$the ${\cal P}_{i}$-cell of $o$ is
intersected by the boundary of some ${\cal P}_{j}$-cell,
$j>i]\leq cb_i/a_i
+\sum_{j=i+1}^\infty c' b_i/b_j
$, with some constants $c,c'$. \Qed

For any $j$, call the cells of ${\cal R}_j '$ that are not
$b_j\times\ldots\times b_j$ cubes {\bf irregular} cells. Note that by \ref
l.seged/, the probability that $o$ is contained in an irregular cell of
some $\rk '$, $k\geq j$, is $\leq c' 2^{-2^j}$. Hence every irregular cell
(in any of the $\rj '$) is contained in some maximal irregular cell of
some $\rk '$ ($k\geq j$), and the probability that this $k$ is greater
than some $\kappa$ is $\leq c 2^{-2^\kappa}$. Let ${\cal I}$ be the set
of maximal irregular cells. We mention that ${\cal I}$ is not necessarily a partition of $\R^d$, but a set of pairwise disjoint subsets of it.
Let $\rj ''$ be the common coarsening of
${\cal R}_j '$,
and ${\cal I}$. 
By \ref l.coarse/, ${\cal R}_j ''$ is still a sequence of coarser and coarser partitions.

%

Let $\rj$ be a refinement of $\rj ''$, to be defined as follows. For each cell $C$ in ${\cal 
I}$, we subdivide $C$
by a $\Z^d$ grid (placed on $C$ in some deterministic way, say 
with origin on an extremal point for some fixed hyperplane). Now replace every cell $C\in {\cal I}\cap {\cal R}_j ''$ by this refinement. The other cells 
of $\rj ''$ (those that are not cells of ${\cal I}$) are unchanged.
Of course we still have a sequence of coarser and coarser partitions, and \ref l.seged/ 
remains valid for the resulting $\rj$:

\procl l.seged2
The probability that $o$ is not in a $b_j\times\ldots\times b_j$ cube of
$\rj $ is $\leq \sum_{k=j}^\infty cb_k/a_k + \sum_{i=j+1}^\infty c_0
b_j/b_i
\leq c2^{-2^i}$. In this case, $o$ is in a cell that is contained in a 
$1\times\ldots\times 1$ cube.
\endprocl

An important fact is
that the partition $\rj$ is completely determined by the elements of the 
$S_i$'s with 
$i=a_j,a_{j+1},\ldots$, and that by \ref l.seged2/ a cell of $\rj$ is 
either a 
$b_j\times\ldots\times b_j$ cube (in which case we call it a {\bf good 
cell}), or a cell of diameter $\leq d^{-1/2}$. The ${\cal R}_j$ satisfy 
the claim of \ref t.main/ (as we show at the end of this section), except 
for that the sizes of the typical cubes 
grow fast (and not just double) as we increase $j$ one by one.  

Given a subcube $H$ of $\R^d$, recall the definition of right-faces and 
$\theta (H)$ (with respect to some given partition of $\Z^d$) from before 
\ref l.lepcsos/. Now we are ready to define the 
final sequence ${\cal Q}_j$ of 
partitions, as in \ref t.main/. The sequence $\{{\cal Q}_j\}$ will be such 
that $\{{\cal 
R}_i\}$ is a subsequence of $\{{\cal Q}_j\}$.

It is enough to define the ``intermediate" partitions between $\ri$ and 
${\cal R}_{i+1}$, for any $i$. Note that the cubic $\ri$-cells do not 
necessarily subdivide $C$ like a cubic grid, as illustrated by the left 
side of 
Figure 1.
For each good cell $C$ of ${\cal R}_{i+1}$ 
(that is, a ${b_{i+1}}\times\ldots\times {b_{i+1}}$ cube), and for 
$\ell=0,1,\ldots, 2^i -1$, consider the 
subdivision $H_\ell$ of $C$ to dyadic cubes of size $b_i 2^\ell$. Now, for 
each 
$K\in H_\ell$, consider $\theta (K)$ with respect to the partition 
$\ri$. Define 
${\cal Q}_i^\ell$ 
as 
the set of $\theta (K)$'s as $K\in H_\ell$. (See the right side of Figure 
1.)

\procl l.valami
The cell $C_o$ of $o$ in $\qi^\ell$ is a $b_i 2^\ell$-pseudocube whenever
$o$ is
in a good cell of ${\cal R}_{i+1}$, unless $\ell =1$ and $o$ is in the
$b_i$-neighborhood of the right boundary of its ${\cal R}_{i+1}$-cell.
This exceptional event has probability $\leq b_i/b_{i+1}=2^{-2^i}$.
\endprocl

\proof
For $\ell=0$ the claim is true even with cubes instead of pseudocubes. The
volume condition for a pseudocube is clear from
the construction and \ref l.isoperi/. To verify the isoperimetry
condition, let $H\in H_\ell$ be the cube that we used to define
$C_o=\theta (H)$,
and apply \ref l.lepcsos/ together with the volume condition.
\Qed

Finally, let $(\qj)$ be the
sequence resulting from the finite sequences $(\qi ^\ell)_\ell$ when we
put them one after the other as $i=1,2,\ldots$

\medskip
\medskip 
\bigskip
\SetLabels
   (.5*-.15) {To the left: The subpartition of a good cell $C$ of ${\cal R}_{i+1}$ by the cells of ${\cal R}_i$. To the right:} \\ 
   (.5*-.25) {the partition $H_\ell$ of $C$ by dyadic cubes (dashed), and 
the $\theta (K)$ (thick lines), $K\in H_\ell$.} \\
  (.5*-.37) {\bf Figure 1.}\\
\endSetLabels
\AffixLabels{\epsfysize=6.3cm
\epsfbox{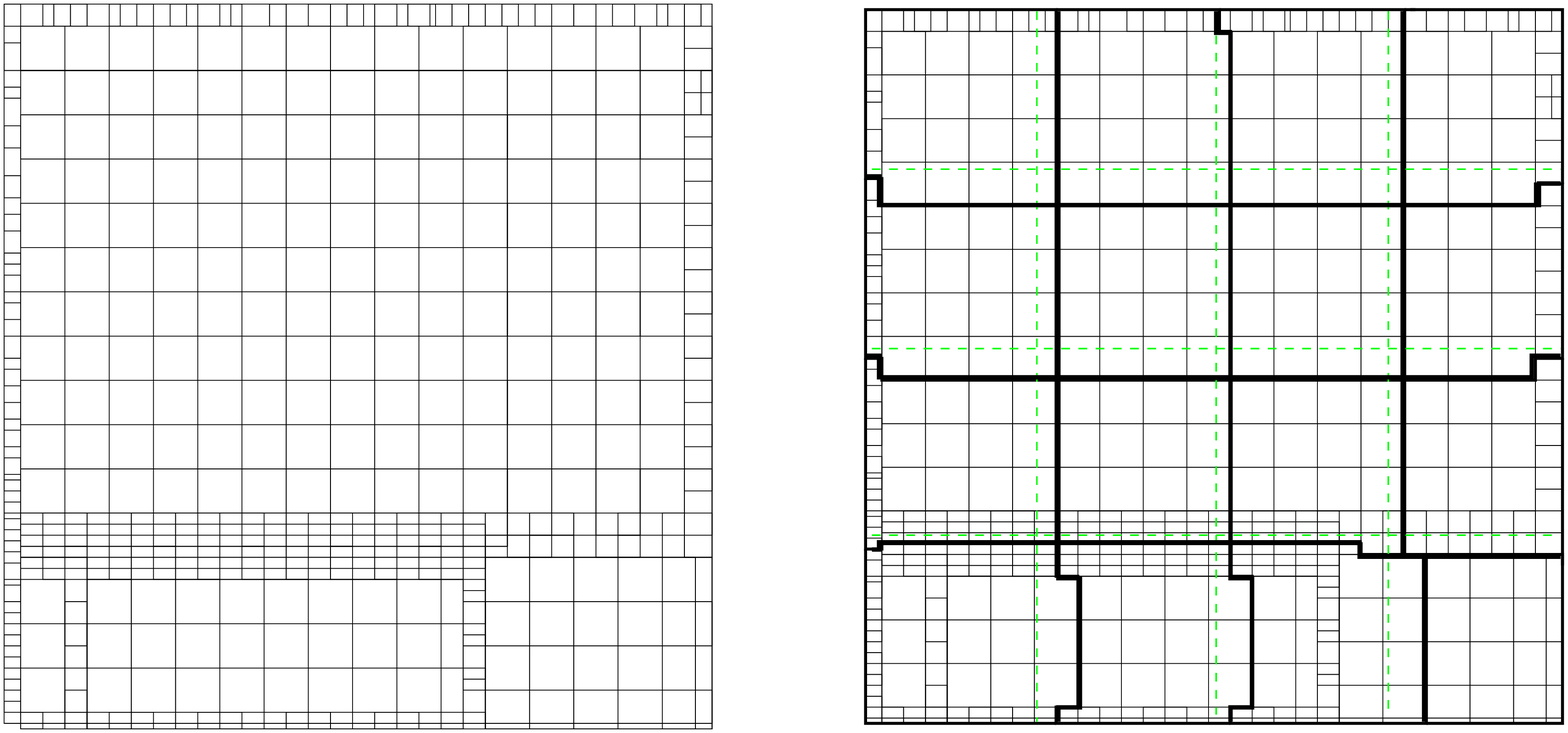}
}

\medskip
\medskip
$$$$
$$$$
$$$$

To finish the proof of \ref t.main/, we shall show that for a fixed 
$o\in\Z^d$ the labels 
of the vertices in the $\qj$-cell of $o$ are i.i.d. conditioned on an 
event of probability tending to 1 rapidly with $j$. 

Before going into that, let us remind ourselves to the construction of 
$\qj$. See Figure 2a and 2b for schematic pictures of the sequentially 
constructed partitions and ${\cal I}$. Note that we changed the scales 
for the figure (and this is the cause of seemingly many noncubic cells, 
which is not the case when one uses the the proper parameters as in our 
construction). 

\procl r.summary {\bf Summary of the construction of the partitions:}
\item{(1)} $S_k$ consists of the vertices that have $k$-bulbs around 
them. 
\item{(2)} $\pk$ is a subpartition of the Voronoi tessallation ${\cal 
V}(S_{a_k})$ on $S_{a_k}$ to
$b_k\times\ldots\times b_k$ cubes (with the exception of a small 
proportion of the cells). 
\item{(3)} $\rk '$ is the 
common refinement of 
the sequence of ${\cal P}_{i}$'s, with index $i$ starting from 
$k$. 
This way most of the cells in $\rk'$ coincided with the cubes in ${\cal 
P}_{k}$; on the other hand the $\rk '$ is a sequence of coarser 
and coarser partitions. 
\item{(4)} ${\cal I}$ is the set of cells $C$ that are not 
$b_k\times\ldots\times b_k$ cubes, with $C\in {\cal R}_k'$, and maximal 
with this property. ${\cal R}_k''$ is the common coarsening of ${\cal I}$ 
and ${\cal R}_k '$.
\item{(5)} We defined $\rk$ from $\rk ''$ by subdividing its non-cubic 
cells (and possibly some others) to small chunks whose 
diameters are uniformly bounded. We still have a sequence of coarser 
and coarser partitions, and most of the cells are still $b_k$ by $b_k$ 
cubes.
\item{(6)} Finally, 
$\qj$ is a sequence that we obtained from 
$\rk$ 
by putting ``intermediate" partitions in the sequence $({\cal R}_k)$ so 
that 
a ``typical cell" is a pseudocube of size always doubling as $j$ increases 
by 1. The subsequence $\rk$ provides the one given in (iii).
\endprocl

Let us point out that the constructions of the partitions ${\cal P}_j, 
{\cal R}_j, {\cal R}_j '$ 
and ${\cal Q}_j$ did not use any information besides that coming from 
$S_j,S_{j+1},\ldots$:

\procl l.pointout
Every partition ${\cal Q}_j$ is a deterministic function of $S_j$.
\endprocl

By \ref l.valami/ we know that $\qj$ satisfies {\it (i)} in \ref t.main/.

To show (ii) in \ref t.main/, fix $j$ and let $i$ be such that 
$b_i<2^j\leq b_{i+1}$, that is, 
$2^i<j\leq 2^{i+1}$. Let 
$B:=B(o,2 b_{i+1})$. Note that the cell of $o$ in $\qj$ 
is 
contained 
in $B$.
Further, $a_i$ is such that 
one was using only elements of $S_{a_i}$ (and hence possibly elements of 
$S_n$ 
with 
$n>a_i$) to define $\ri$ (and thus $\qj$).
Call a vertex $x$ in $\Z^d$ a {\bf giant}, if 
it is an $a_i$-giant, as defined in \ref l.giant/ 
($t=a_i$).
Note that by definition, the existence of giants is independent of the 
configuration within $B$.

Our key observation is that if
$\Z^d$ contains no giants, then no element $x$ of $S_{a_i}$ can be so
close to
$o$ as that the largest bulb around $x$
intersects $B$. That is, if there are no giants, one can tell the elements
of $S_{a_i},S_{a_{i+1}},\ldots$ without looking into $B$.
Hence the
configuration outside $B$ determines the cell $C$ of $o$ in ${\cal
R}_{i+1}$, together with the subpartition of $C$ by $\ri$ - and these two
determine the cell of $o$ in $\qj$. (Here we are using \ref l.pointout/.)
Then, conditioned on this event (no giants), the vertices in the cell of
$o$ in $\qj$ have
i.i.d. labels, since the cell of $o$ is contained in $B$.
Now, by \ref l.giant/, there are no giants with probability $\geq
1-Cc^{-a_i^{d}}$. Furthermore, the probability that the cell of $o$ in
$\qj$ is not a $2^j\times\ldots\times 2^j$ pseudocube is at most the
probability that it is not in a $b_i\times\ldots\times b_i$ cube of
$\ri$ plus the probability that it is at the right boundary of its
$\ri$-cell, as in \ref l.valami/. By \ref l.seged2/ and \ref l.valami/, this 
is bounded by
$c2^{-2^i}+c'b_i^{-1}\leq c''2^{-2^i}\leq c''2^{-2^{j-1}}$.

Define $A_j$ to be the event that
there are no giants, and the cell of $o$
is a
pseudocube. We have just seen that $P[A_j]\geq 1-c2^{-2^{j-1}}$, and this
finishes the proof of (i) and (ii) in  \ref t.main/. Part (iii) follows by
setting ${\cal Q}_{\alpha(k)}:={\cal R}_k$.

\SetLabels
   (.5*-.15) {\bf Figure 2a.} Construction of the sequence of partitions\\
   (.52*-.18) (on a scale altered from the real one)\\
   (.2*.55) Four points of $S_k$, their\\ 
   (.2*.52) neighborhood, and ${\cal 
    V}(S_k)$\\
   (.2*-.05) ${\cal R}_k '$ (thicker lines indicate boundaries\\
   (.2*-.08) of ${\cal 
    P}_j$, $j=k+1,\ldots$)\\
   (.78*-.05) ${\cal I}$ (the grey regions are not \\
   (.78*-.08) covered by any set of 
   ${\cal I}$)\\
   (.78*.55) ${\cal P}_k$\\
  \endSetLabels
\AffixLabels{\epsfysize 16cm
\epsfbox{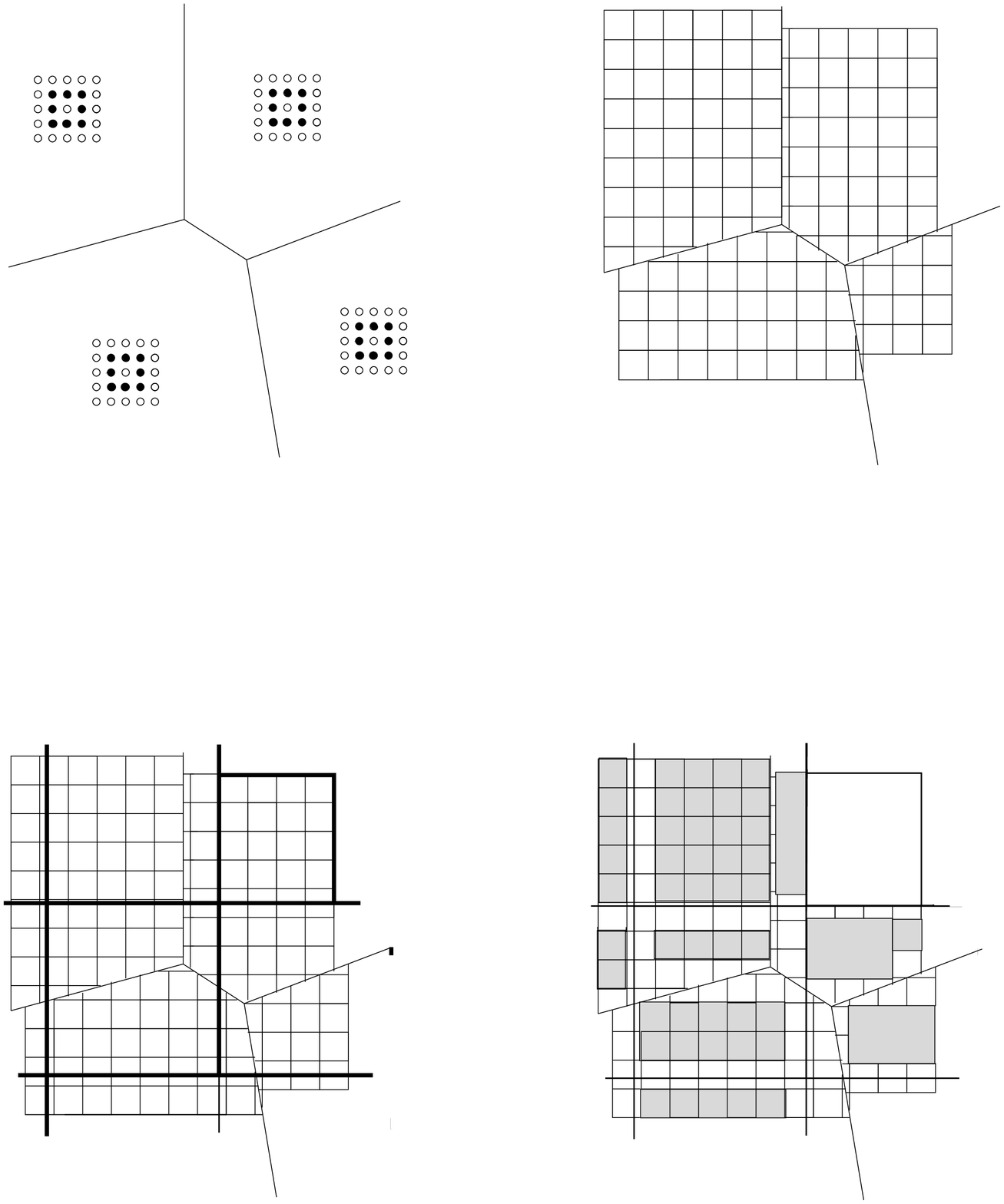}}
$$$$
$$$$

{\hfill
\SetLabels
   (.5*-.09) ${\cal R}_k$, resulting from ${\cal R}_k ''$ (the finest common coarsening of ${\cal I}$ and ${\cal R}_k '$) \\  
   (.5*-.15) by subdividing irregular cells to cells of small diameters \\
   (.5*-.25) {\bf Figure 2b.}\\
  \endSetLabels
\AffixLabels{\epsfysize 8cm
\epsfbox{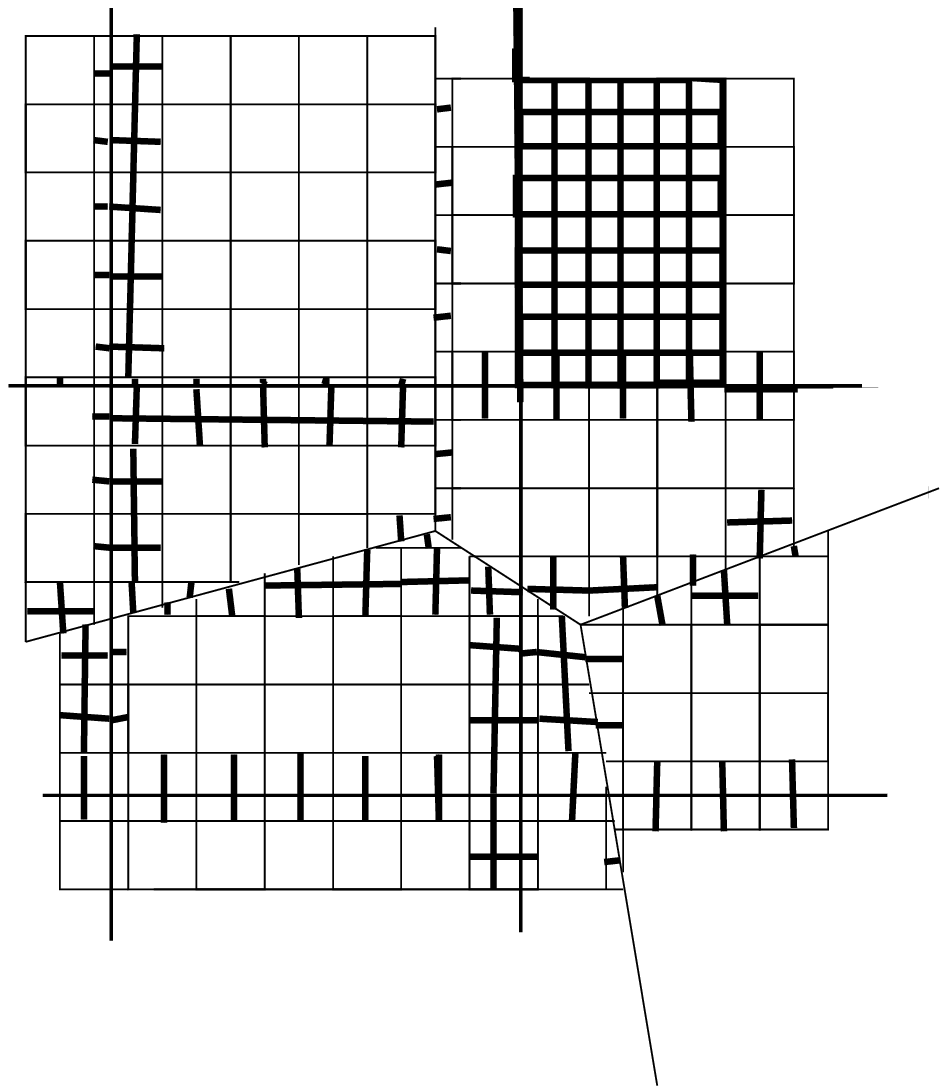}
}\hfill}
\bigskip
$$$$
$$$$

\bsection{The matching rule}{s.matching}

Recall that an $n$-{\it pseudocube} or a {\it pseudocube of size} $n$ in 
$\Z^d$ is a subset $C$ that contains some $n/2\times\ldots\times n/2$ 
cube and is contained in some $2n\times\ldots\times 2n$ cube, and 
satisfies the isoperimetry condition. Call the elements of $\cup_i {\cal 
Q}_i$ {\bf dyadic pseudocubes}. The reason for the name is that by \ref 
t.main/, the  pseudocube of ${\cal Q}_i$ is about twice the size of 
the pseudocube in ${\cal Q}_{i-1}$, and these partitions are 
coarser and coarser, so most of the pseudocubes in ${\cal Q}_i$ are 
subdivided by pseudocubes of ${\cal Q}_{i-1}$ in a dyadic pattern.

For some subset $S$ of $Z^d$, denote by $\y (S)$ the set of yellow
elements of $S$, and by $\b (S)$ the set of blue elements of $S$.

Given some $C\subset\Z^d$, with
Bernoulli($1/2$) labels on it, let the {\bf surplus} of $C$ be 
$|\b (C)-\y (C)|$. Denote this quantity by $\sur (C)$.

Note the distinction between {\it subgraphs} of a graph $G$ and {\it
graphs on the vertex set of} $G$. Also, there is a bit of ambiguity about
the use of the word {\it edge} : sometimes it refers to edges of $\Z^d$,
and sometimes to pairs in the matching, but it is always clear from the
context.

From now on {\it fix function} $f:\R^+\mapsto \R^+$ to be 
$f(x)=x^{{d-1-\epsilon/2\over d}}$.
Say that a pseudocube $C$ is {\bf bad}, if $\sur (C)>f(|C|)$.

Before presenting the matching rule, let us prove a few simple lemmas. 

\def\min{{\rm min}}
\def\dist{{\rm dist}}
\def\diam{{\rm diam}}

\procl l.bound
Let $\Delta\subset\R^d$ be connected, and $U\subset \R^d$ with 
$\min_{x,y\in U} |x-y|\geq s$. Then if $\diam (\Delta) < s$, then $|U\cap 
\Delta|\leq 1$. Otherwise 
$$|\Delta\cap U|\leq c|\partial \Delta|\diam (\Delta) s^{-d}.$$
\endprocl

\proof 
If $\diam (\Delta)\leq s$, then the statement trivially holds by the 
assumption on $U$. So suppose $\diam (\Delta) > s$. Let 
$\Delta_s$ be the $s$-neighborhood (``fattening") of $\Delta$. We have 
$|\Delta_s|\leq c |\partial 
\Delta|\diam (\Delta)$. Since the balls of radius $s/2$ around points 
of $\Delta\cap U$ are disjoint, and they are all contained in $\Delta_s$, 
we have $|\Delta\cap U|\leq c'|\Delta_s|s^{-d}\leq c |\partial 
\Delta|\diam 
(\Delta)s^{-d}$.
\Qed 

\procl l.sparse Let $S\subset \R^d$ be a cube, $S'\subset S$, $|S'|\geq 
|S|/2$. Then there exists an $S''\subset S'$ with 
$|S''|=f(|S|)=|S|^{{d-1-\epsilon/2\over d}}/2$ such that $\dist 
(x,y)\geq c |S|^{{1+\epsilon/2 \over d^2}}$ for every $x,y\in S''$, with $c$ depending only on 
$d$. \endprocl

\proof
Pick
elements of $S'$ for $S''$ one by one, always removing points of the 
$c|S|^{{1+\epsilon/2 \over d^2}}$-neighborhood of the chosen point from 
$S$ 
(and picking next elements 
from what remains).  \Qed

Fix a $k$, to be determined later. We will choose it large 
enough, and a power of 2 for technical 
convenience. Note that in what follows, $c$ and $c'$ will always denote 
constants that {\it do not dependent on} $k$. 
Define $a (0):=\log k$, and $b(1):=1$.

{\bf Choice of} $a(i)$ {\bf and} $b(i)$: 
Fix $b(i)>2^{id}$,
$i=2,3,\ldots$, to be a subsequence of $\alpha(i)$ 
(where
$\alpha (i)$ is as
defined in \ref
t.main/).
Fix an increasing subsequence $a(i)$ of $\alpha(i)$, $i=1,2,\ldots$, such 
that
$a(i)>\max \{b(i)^{2d/\epsilon},\exp (\exp (i^{d-2}))  \}$, and such 
that 
$\sum_{\nu=a(i)/d}^\infty (k2^\nu)^{-\epsilon/2}\leq {1\over 2} k^{{1\over 
d-1}}b(i+1)^{-d}$.

The next lemma for later use is of elementary geometry. The claim 
is intuitively clear, but not trivial to prove. We do not need the 
specific value of the constant, but $c_d=2^{d+2}$ works. 
Recall, that for any 
subset 
$S$ of $Z^d$, we denote by $\partial S$ 
the 
set of vertices with degree 
$<2d$ in $S$. In the next lemma, if $E$ is some set of edges and $V$ is 
some set of vertices, 
let $E_V$ stand for the subset of edges $E$ incident to some element 
of 
$V$.

\procl l.surface
There exists a constant $c_d$, depending only on $d$, such that the 
following holds.
Let $\Gamma$ be some subset of the edges of an $m\times\ldots\times m$ 
cube $K$, and let $\Delta_1,\ldots,\Delta_h$ be the connected components 
of 
$K\setminus \Gamma$. 
Then for all but at most one of the $j$'s we have
$$|\partial \Delta_j|\leq c_d |\Gamma_{\partial\Delta_j}|.$$
Consequently,
there is a $\Delta_i$ such that
$$\sum_{j\not=i}|\partial\Delta_j|\leq
2c_d |\Gamma|.$$\endprocl

\proof
Let $i$ be such that $|\partial\Delta_i\cap\partial K|$ is maximal (in 
case of ambiguity, decide arbitrarily). Take any $j\not= i$.
For $x\in\partial K$ and $z\in K$, define $P(x,z)$ as a path from $x$ to 
$z$
that makes as few turns as possible (otherwise its choice is arbitrary).
Regard $P(x,y)$ as some element of the vector space 
$\R^{\overrightarrow{E}(K)}$, where ${\overrightarrow{E}(K)}$ is some 
arbitrary fixed orientation of the edges of $K$. For each element $(x,y)$ 
of $(\partial\Delta_j\cap\partial K)\times(\partial 
K\setminus \partial\Delta_j)$, let
$P_{x,y}\in\R^{\overrightarrow{E}(K)}$ be ${1\over |K|}\sum_{z\in K} 
(P(x,z)-P(z,y))$. 
Finally, define a flow $\mu$ 
as 
$\sum_{x,y} {1\over |\partial
K\setminus \partial\Delta_j|}P_{x,y}$, where the sum is over  
$(x,y)\in (\partial\Delta_j\cap\partial
K)\times(\partial
K\setminus \partial\Delta_j)$. It is not hard to see that the 
flow through any edge is $\leq {(c_d-1)\over 2}(1+ 
{|\partial\Delta_j\cap\partial
K|\over  |\partial
K\setminus \partial\Delta_j|})\leq (c_d-1)$ with some constant $c_d$, 
where the second inequality follows from $j\not= i$. (A 
proof 
for the 2-dimensional 
case, can be found  
in \ref b.A/, Lemma 9. The only additional thing needed, is an upper bound 
of order $|K|$
for the number of paths $P(x,z)$, $(x,z)\in \partial K\times K$ containing 
an arbitrary edge $e$.
Now, 
if 
$e=\{(v_1,\ldots,v_\nu,\ldots,v_d),(v_1,\ldots,v_\nu+1,\ldots,v_d)\}$,
then
one of $x$ and $z$ has $j$'th coordinate equal to $v_j$ for every
$j\not=\nu$, which gives an order $|K|^{1+1/d}$ choices for them, but
also $x$ has to be on the boundary, which makes the number of choices be
of order $|K|$. ) 
The total 
strength of the flow is 
$|\partial\Delta_j\cap\partial
K|$. Since $ \Gamma_{\partial\Delta_j}$ is a cutset between the sources 
and the sinks of the flow, the total amount flowing through it is equal 
to the strength. Putting this together with the bound on the amount 
flowing through an edge, we obtain:
$$|\partial\Delta_j\cap\partial
K|\leq (c_d-1)| \Gamma_{\partial\Delta_j}|.$$
Hence $|\partial\Delta_j|\leq |(\partial\Delta_j\cap\partial
K)\cup \Gamma_{\partial\Delta_j}|\leq 
c_d |\Gamma_{\partial\Delta_j}|$.\Qed

The following lemma of several later uses is Chernoff's bound about 
independent Bernoulli sums, in the language of surpluses.
                                                                                
\procl l.chernoff
The probability that a fixed pseudocube $K$ of size $k2^i$
has surplus $\geq f ((k2^i)^d)$ is $<\exp (-f((k2^{i-1})^d)^2/2(k2^{i+1})^d)$.
\endprocl

\proof
Chenoff's bound provides an estimate $<\exp ({-f(|K|)^2\over 2|K|})$ for 
the 
probability. Then use the fact $(k2^{i-1})^d\leq |K|\leq (k2^{i+1})^d)$. \Qed



{\bf Definition of the boundary of a partition:}
Given some partition ${\cal P}$ of the vertices of some graph in $\Z^d$, {\it we 
denote 
by 
$\partial {\cal P}$ the subgraph induced by the union of the (inner) vertex boundaries of the 
cells of 
${\cal 
P}$}.

{\bf Definition of} ${\cal P}$:
Consider ${\cal Q}_i$ from \ref t.main/.
For a $C\in {\cal Q}_i$, let the cells of ${\cal Q}_{i-1}$ in $C$ be 
called the {\bf bricks} of $C$.
Let the cells of 
${\cal Q}_{\log k}={\cal Q}_{a(0)}$ be called {\bf basic pseudocubes}. 
For $j\geq \log k$, say that a cell $C\in{\cal Q}_{j}$ is {\bf 
bad}, if 
$\sur (C)\geq f(|C|)$, or if $C$ is not a $2^j$-pseudocube. This 
definition implicitly relies on $j$ and ${\cal Q}_{j}$; however, 
one can extend the definition of ``bad" to any element $C$ of $\cup 
{\cal 
Q}_j$ for the following reason. Since the ${\cal Q}_j$ are coarser and 
coarser partitions, one can trace back for each $C\in\cup_{j=\log k}^\infty 
{\cal Q}_j$ the smallest $j$ such that $C\in {\cal Q}_j$, and say that 
$C$ is bad, if it is bad in ${\cal Q}_j$.

Call $C$ 
{\bf ripe}, 
if $C$ is not bad, but one of its bricks is 
bad, 
and further, $C$ is maximal among pseudocubes of $\cup_j \qj$ with this 
property with respect 
to 
inclusion. That in fact there is a maximal such $C$ containing $o$ for 
every $o$, follows from the next lemma.

\procl l.chernoff2
The probability that ${\cal Q}_{i+\log k} (o)$ is not bad, but one of its 
bricks is 
bad is at most $c \exp (-c'(k2^i)^{d-2-\epsilon})$.\endprocl

\proof
For $x\in C$, the probability that ${\cal Q}_{i-1+\log k}(x)$ is bad is 
bounded by $c \exp (-c'(k2^i)^{d-2-\epsilon})$. This is a consequence of 
\ref 
l.chernoff/ (in case ${\cal Q}_{i-1+\log k}(x)$ is bad because it has too 
large
surplus), and \ref t.main/ (in case ${\cal Q}_{i-1+\log k}(x)$ is bad 
because 
it is not a pseudocube). Hence
$$\P[{\cal Q}_{i+\log k} (o) \hbox{ is not bad but some brick of it is 
bad}]\leq $$
$$\leq\P[{\cal 
Q}_{i-1+\log k}(x)\hbox{ is bad for some }x\in B_{k2^{i+2}}] \leq c \exp 
(-c'(k2^i)^{d-2-\epsilon})$$
if $k$ is large enough, for every $i\in\Z^+$.
\Qed

A bound of the same 
magnitude 
holds for the probability of being in a ripe pseudocube as for the 
probability 
of being in a bad pseudocube of similar size, as stated in \ref 
l.chernoffnew/. 

Call a pseudocube $C$ an {\bf elementary pseudocube} if it is ripe or it 
is a basic pseudocube that is not contained in any ripe pseudocube.
It's a consequence of the definitions that elementary pseudocubes in 
$\cup_j 
{\cal Q}_j$ 
partition $\Z^d$; call this partition ${\cal P}$. 
An important property
of ${\cal P}$ is that none of its cells $C\in{\cal Q}_j$ is bad in ${\cal 
Q}_j$.
\ref l.chernoffnew/ 
gives an exponential bound on the tail probability of the diameter of 
the ${\cal P}$-cell of a point. In fact, the reason we defined 
pseudocubes, and the partitions of always doubling approximate cell sizes 
in \ref t.main/ (instead of just taking the more convenient sequence 
${\cal Q}_{\alpha (k)}$), is to have this control of the tail coming from 
\ref l.chernoff/ (which is 
inherited from the tail for bad cubes).

\procl l.chernoffnew
The probability that a vertex is in a ripe pseudocube of size $k2^i$ is
at most $c \exp (-c'(k2^i)^{d-2-\epsilon})$. Hence, the probability that 
the
${\cal P}$-cell of a point has size $k2^i$ has the same bound.
\endprocl

\proof
By \ref l.chernoff/ and the fact that a ripe pseudocube of size $k2^i$
contains a bad pseudocube of size $k2^{i-1}$.\Qed

{\bf Definition of} ${\cal P}_i$:
If we find a maximal matching within each ${\cal P}$-cell, the set of 
unmatched points 
is ``relatively small". We will then match these points, using a 
sequence of coarser and coarser partitions.  By definition, ${\cal P}$ 
is such that for any cell $C$ of ${\cal Q}_i$, $i\geq \log k$, $C$ is either contained in 
some cell of ${\cal P}$, or it is a union of some cells of ${\cal P}$.
Define ${\cal P}_i$ as the 
common coarsening for ${\cal Q}_i$ and ${\cal P}$, $i\geq a(0)$. (In 
particular, ${\cal P}_0={\cal P}$.) Since ${\cal 
P}$ is a refinement of ${\cal P}_i$ for any $i$, every cell of ${\cal 
P}_i$ contains some elementary cell. In 
particular, we achieve the following, which was the goal of the last few paragraphs:

\procl l.no 
If $i\geq \log k$ then 
no cell of ${\cal P}_i$ is bad. That is, every cell $C\in {\cal P}_i$ is a 
$2^{i'}$-pseudocube with some $i'\geq i$, and $\sur (C)\leq f(|C|)$. 
\endprocl

Call a cell $C\in {\cal 
P}_j$ {\bf cubic}, if it is a $2^j$-cube, otherwise it is {\bf non-cubic}. 
If $j=\alpha (i)$ for some $i$, then every cell of ${\cal Q}_j$ that is 
not bad is cubic, hence ${\cal P}_j (o)$ can be non-cubic only if it 
is a cell in ${\cal P}\setminus {\cal Q}_j$, and hence ${\cal P}_j (o)$ 
has 
size $\geq 2^j$. 
This probability is bounded by  
\ref l.chernoffnew/. These two give us that with probability exponential 
in $2^{j(d-2-\epsilon)}$, a point is in a $2^j$-cubic cell of $C\in{\cal
P}_{j}={\cal P}_{\alpha (i)}$.

For each $i$ we will define an invariant monochromatic subset $U_i$ of 
$Z^d$ and a 
matching ${\cal M}_i$ in such a way that $|U_i\cap C|=sur(C)$ for every 
$C\in 
{\cal P}_{a(i)}$, and every vertex of $C\setminus U_i$ is matched by 
${\cal M}_i$.

For every $2^{a(i)}\times\ldots\times 2^{a(i)}$ cube $S$ in ${\cal 
P}_{a(i)}$, if the majority of points in $S$ is blue, choose $S'\subset 
S$ to be the 
set of blue vertices, otherwise let $S'$ be the set of yellow vertices. 
Choose $S''$ from $S'$ as in \ref l.sparse/. Finally, pick an arbitrary subset of 
cardinality $\sur (S)$ from $S''$; call this $U_i (S)$. Now let $U_i:=\cup 
U_i (S)$.
(Of course 
the 
above 
choices can be made cell by cell in some predetermined way, that is the 
same for every translate of the cell, to make the 
resulting $U_i$ invariant.)
Note that 
$$\min_{x,y\in U_i\cap S} |x-y|\geq c|S|^{{1+\epsilon 
/2\over d^2}}=2^{{a(i)(1+\epsilon
/2)\over d}},  \label e.kiemel
$$
by our choice coming from \ref l.sparse/.

Then we have, by \ref l.bound/, the following:

\procl l.combined
For any $S$ cubic cell of size $2^{a(i)}$ of ${\cal P}_{a(i)}$, and 
connected subset $\Delta$ of $S$, 
either $|\Delta\cap U_i|\leq 1$, or
one has 
$$|\Delta\cap U_i|\leq c|\partial\Delta |2^{a(i)}(\min_{x,y\in 
U_i}\dist(x,y))^{-d}\leq c |\partial \Delta|2^{-\epsilon a(i)/2}.$$

\endprocl

\bigskip
\SetLabels
  (.25*-.12) $K \in{\cal P}_{a(i)}$\\ 
  \L(.02*-.22) White cells represent cells of ${\cal P}$.\\
  \L(.02*-.29) Thick lines show cell boundaries of ${\cal P}_{b(i)}$.\\
   (.5*-.45) {\bf Figure 3.} The contraction $\tilde \beta$ from $K$ to $\tilde K$. (Scale altered from real. Pseudocubes of \\
   \L(.1*-.53) ${\cal P}$ and ${\cal P}_{b(i)}$ are represented by cubes.)\\
  \L(.75*-.12) ${\tilde K}$\\
  \L(.54*-.22) The multiplicity of each edge is the number\\
  \L(.54*-.29) of edges between two small cubes on the left.\\
\endSetLabels
\AffixLabels{\epsfysize=6.3cm
\epsfbox{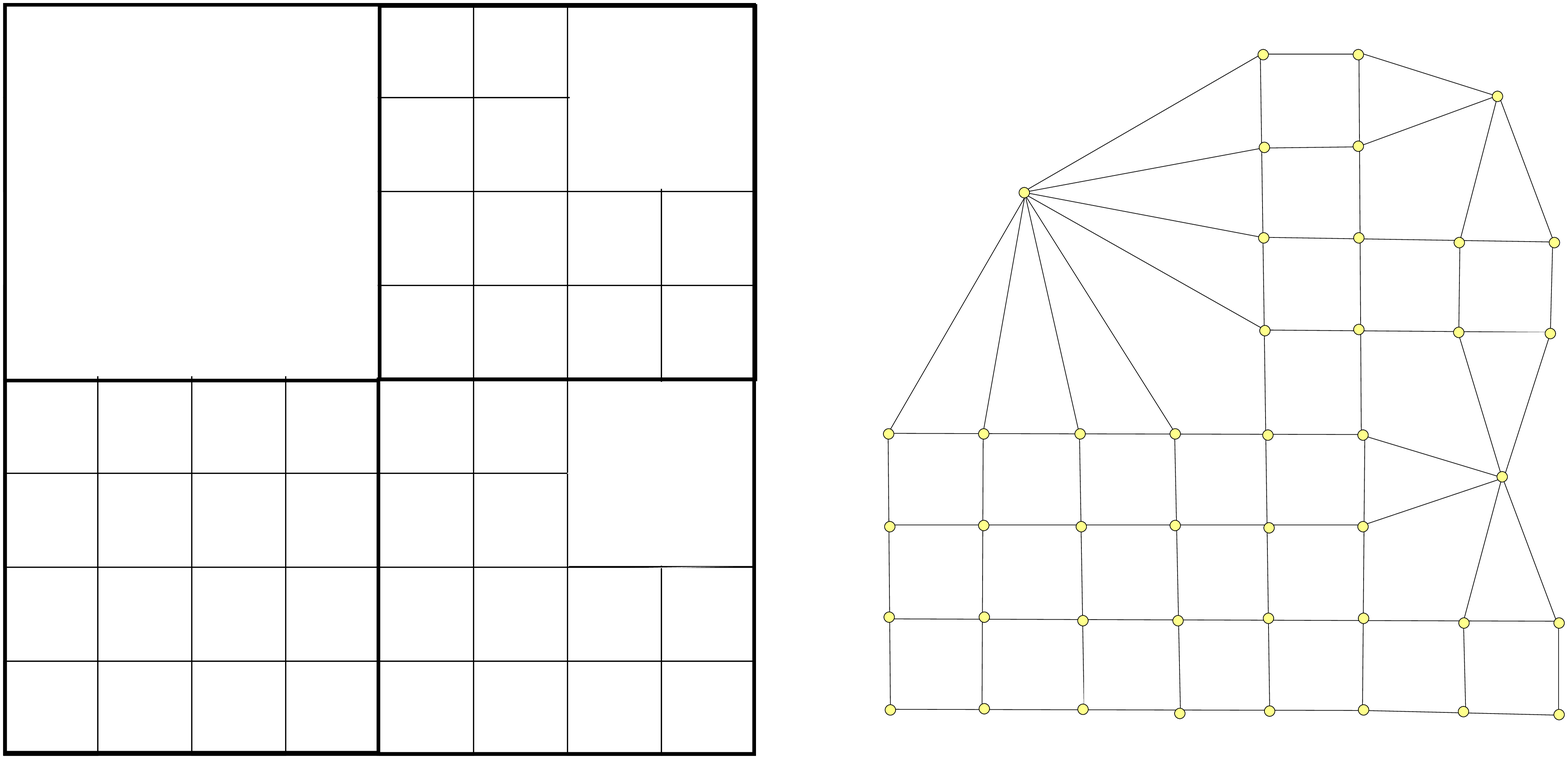}
}

\medskip
\medskip
$$$$
$$$$
$$$$
$$$$

Now we are ready to present the {\bf matching algorithm}.

{\bf Step 1.
}
For each elementary pseudocube $C$ (i.e. cell of ${\cal P}$) match
blue points with
yellow points as long as it is possible.

Take the union of
these matchings over all elementary pseudocubes $C$, and call the 
resulting
matching
${\cal M}_0$. Note that in each elementary pseudocube $C$ the number of 
points
not
matched by ${\cal M}_0$ is equal to $\sur (C)\leq f(|C|)$ (using \ref 
l.no/ applied to ${\cal P}_{a(0)}$, i.e. 
the fact 
that
an elementary pseudocube is not bad). Hence, at least 
$|C|/3$
vertices of each color in every
elementary pseudocube $|C|$ are matched by ${\cal M}_0$.
This, and the bound on the number of unmatched points will have importance when we will define pairwise disjoint augmenting paths.
It is also clear that ${\cal M}_0$ is invariantly defined.

The next lemma is straightforward from \ref l.chernoffnew/ and the 
definition of ${\cal M}_0$.

\procl l.M0 
The probability that a vertex, given that it is matched by ${\cal M}_0$, 
has its 
pair at distance $\geq r$ from it, is 
$\leq c \exp (-c' f(r^d)^2/r^d)
\leq c \exp (-c'r^{d-2-\epsilon}).$
\endprocl

Denote by $U_0$ the set of points that are not matched by ${\cal M}_0$. 

For sets $A$ and $B$, the {\it multiset} union of them will be denoted by 
$A\dot\cup B$; that is, $A\dot\cup B$ is the union of $A$ and $B$, with 
elements in $A\cap B$ having multiplicity 2. If $A$ is a multiset and $B$ 
is a set, then we define $A\cap B$ to be the multiset such that every 
$x$ that is in $B$ will 
have the same multiplicity in $A\cap B$ as in $A$.

In {\bf Step 2}, we proceed in countably many stages. In {\bf Stage i}, as 
$i=1,2,\ldots$, we define a rematching procedure, that will match every 
point of $\Z^d\setminus U_i$ to some other point. In particular, it 
matches the points of $U_{i-1}\setminus U_i$ (which were unmatched 
before). The idea is that in 
Stage $i$ we redefine only edges that have an endpoint in a 
set of ${\cal P}$-cells (which is ``sparse'' of density about $1/b(i)$). The method ensures 
that
the tail of the 
matching 
remains of 
the same magnitude as for ${\cal M}_0$ (provided by \ref l.M0/). The 
rematching will 
use augmenting paths. When $i>1$, in 
Stage $i$, the scarcity of $U_{i-1}$ 
(which is a result of our choice for the sequence $a(i)$) makes the 
rematching simpler. However, in Stage 1, the set $U_0$ of points to be 
matched is coming from Step 1, and this fact is responsible for 
more difficulties and the sharp
role of isoperimetry in this case. 

{\bf Definition of} $N(\tilde K)$: 
For each cubic cell $K$ in ${\cal P}_{a(i)}$ we will do the following. 
{\it Fix} $i$ and a $K$ in ${\cal P}_{a(i)}$ for the rest of this 
section. 
Define $\tilde 
K$ by contracting every cell of ${\cal P}$ to a vertex, that is, 
contracting elementary pseudocubes of $K$. We do not erase any edge after 
the contraction, multiple edges are allowed.
Let $\tilde\beta$ be the contraction mapping from $K$ to $\tilde K$. 
(See Figure 3.)


Then 
define $N(\tilde K)$ to be a network on a graph $(V,E)$, with capacities 
on edges. 
Here 
$V$ is the union of $V(\tilde 
K)$ and two extra vertices $B$ and $Y$, and $E$ is the union of  
$E(\tilde K)$ and all edges of the form $\{B,x\}$ or $\{Y,x\}$, $x\in 
V(\tilde
K)$. Define capacities as follows. Let 
every edge incident to some vertex in 
$\tilde K\setminus 
\partial \tilde\beta  ({\cal P}_{b(i)})$
have capacity 0. 
Recall that $\partial\tilde\beta ({\cal P}_{b(i)})$ denotes the union of 
the vertex boundaries of the cells in the push-forward partition $\tilde\beta ({\cal 
P}_{b(i)})$. Let each edge induced by
$V(\tilde K)\cap 
\partial 
\tilde\beta  ({\cal P}_{b(i)})$ have capacity $k^{{d\over d-1}}/12d'$, 
where $2d> d'\geq d$ is chosen so that it makes this number an integer 
(and at least 1, using that $k$ is large enough). Note that for any 
$x\in\tilde K$, $\sum_{x\in e} \capa (e)$ is at most half the 
number of points in $\tilde\beta^{-1} (\tilde K)$ that are matched by 
${\cal M}_0$, by the definition of ${\cal M}_0$ and $\tilde K$, and the  
properties defining a pseudocube. The same will be true if we replace ${\cal 
M}_0$ by ${\cal M}_{i-1}$, because the endpoints of ${\cal M}_{i-1}$ 
will always contain the endpoints of any previous matching.

Finally, for every vertex $x$ 
in $(U_i\dot\cup U_{i-1})\cap 
K$, choose a ``representative" $\bar x\in  V(\tilde K)\cap \partial
\tilde\beta  ({\cal P}_{b(i)})$ such that $x$ and $\tilde\beta^{-1}\bar 
x$ are in the 
same ${\cal P}_{b(i)}$-cell. 
For each such $\bar x$, if $x$ is blue and $x\in U_{i-1}$ or $x$ is 
yellow and $x\in U_i$, add an edge of capacity 1 
between $B$ 
and $\bar x$, otherwise add an edge of capacity 1
between $Y$
and $\bar x$. Let the other edges on $B$ and $Y$ have capacities 0. 
Note that the total capacity of the edges on $B$ is equal to the total 
capacity of the edges on $Y$, by the choice of $U_i$
(because $\b (U_i\cap K)-\y (U_i\cap K)=\b (U_{i-1}\cap K)-\y(U_{i-1}\cap 
K)$).

What we are really interested in is the set of vertices incident to some 
edge of nonzero capacity in 
$\partial\tilde\beta ({\cal P}_{b(i)})\cup\{B,Y\}$, and hence the network 
induced by them. Having the bigger 
network here, with many edges of zero capacity, has the advantage that it 
is easily related to $K$. This technical convenience will make it easy to 
use
some of our geometric lemmas.

For an $X\subset \tilde K$, denote by $N_{B,X}$ and $N_{Y,X}$ the 
subnetwork of $N(\tilde K)$ induced by $B\cup X$ and $Y\cup X$ 
respectively.

We will form a set $E(i,K)$ of pairs from the elements of $U_{i-1}\dot\cup 
U_i$ in $K$, and 
a set of paths $\tilde{\cal P}(i,K):=\{\tilde P(x,y)\, :\,
\{x,y\}\in E(i,K)\}$, in 
such a way 
that:

\item{(i)} elements of $U_i$ are paired with elements of the same color 
in 
$U_{i-1}$, and all the elements of $U_{i-1}$ that are not paired this way, 
are paired with an element of the opposite color in $U_{i-1}$. Call the 
set 
of such pairs $E(i,K)$.
For any pair $\{x,y\}\in E(i,K)$, 
$\tilde 
P(x,y)$ is a path in
$\tilde K\cap\partial\tilde\beta ({\cal P}_{b(i)})$ between $\bar x$ and
$\bar y$ (where $\bar v$ is defined from $v$ as in the definition of 
$N(\tilde K)$).

\item{(ii)} Every vertex $\tilde x\in \tilde K$ is
contained in at most $|\tilde\beta^{-1}(\tilde x)|/3$ elements of
$\cup_{i,K}\tilde{\cal P}$.

\procl p.NAGY
If the maximal flow on $N(\tilde K)$ from $B$ to $Y$ has  
$\sigma={|U_{i-1}\cap K|+|U_i\cap K|\over 2}$ strength, then there exists 
a set of paths $\tilde{\cal P}(i,K)$ 
(with a set $E(i,K)$ for the pairs of endpoints), that satisfy (i) and 
(ii). 

If such $\tilde{\cal P}(i,K)$ 
exists, then there is a set ${\cal P}(i,K)$ of pairwise vertex-disjoint 
paths on the vertices of $K$ such that for every $\tilde P\in\tilde{\cal 
P}(i,K)$, $\tilde P=(x_1,\ldots, x_m)$, there is a $P\in {\cal P}(i,K)$, 
$P=(v_1,\ldots, v_{m'})$, such that $v_1\in U_{i-1}$, and:

\item{(I)} $m'=2m+1$ if and only if $v_{m'}\in U_i$;
$m'=2m+2$ if and only if $v_{m'}\in U_{i-1}$;

\item{(II)} $v_1$ and $v_{m'}$ are such that $\tilde\beta\bar v_1=x_1$ and $\tilde\beta\bar 
v_{m'}=x_m$ ;

\item{(III)} for $\nu=1,\ldots, m$, one has $v_{2\nu}\in \tilde\beta^{-1} 
(x_\nu)$ and further, $v_{2\nu}$ and $v_{2\nu +1}$ are matched by ${\cal 
M}_{i-1}$. Moreover, the vertex colors along $P$ are alternating.
\endprocl

Figure 4 illustrates the connection between the elements $\tilde{\cal 
P}(i,K)$ and ${\cal P} (i,K)$.

\proof
If there exists an 
admissible flow from $B$ to $Y$ of strength 
${|U_{i-1}\cap K|+|U_i\cap K|\over 2}=\sum_x \capa (\{B,x\})=\sum_x \capa 
(\{Y,x\})=\sigma$, then there is also an integer valued flow, 
since the constraints on every edge are integers. Such a flow
can be decomposed as a sum of paths. Delete the first and last edge (the 
ones incident to $B$ and $Y$) from each of these paths, and define the 
set of resulting paths to be $\tilde{\cal P} (i,K)$. By the 
definition of 
$N(\tilde K)$, this shows the first assertion, (i) and (ii). 
 
We can take a preimage in $K$ by $\tilde \beta$ for each 
element $\tilde P(x,y)$ in $\tilde{\cal P}(i,K)$, increase them by attaching one point of $U_{i-1}\cup U_i$ to each endpoint of the preimage,
to find a 
set of pairwise 
disjoint alternating augmenting paths $P(x,y)$ that satisfy (I)-(III). 

The connection between $P(x,y)$ and $\tilde P(x,y)$ is simply the 
following. The 
two endpoints of $P(x,y)$ are $x$ and $y$. 
The $2j$'th point of $P(x,y)$
is going to be a vertex from $\tilde\beta ^{-1}(x_j)$ that is covered
by ${\cal M}_{i-1}$, and chosen to be of a color different from $v_1$. 
Then the $2j+1$'st vertex will be the pair of $v$ by
${\cal M}_{i-1}$. The choices of the $v_{2j}$ are otherwise arbitrary, 
except for that ${\cal P}(i,K):=\{P(x,y)\, :\, \{x,y\}\in E(i,K)\}$ has to 
consist of pairwise disjoint paths.

We
can indeed choose the preimages
$P(x,y)$ to be pairwise
disjoint and fulfill (III), by the choice of the
capacity constraints in $V(\tilde K)$: every $x\in V(\tilde K)$ is crossed
by at most as many paths of $\tilde{\cal P}$, as the number of edges of
$M_{i-1}$ with a yellow (respectively: blue) endpoint in
$\tilde\beta^{-1}(x)$, and so there is a  
choice when each of these edges is present in at most one of the
augmenting paths.
\Qed

\medskip
\medskip 
\bigskip
\SetLabels
   (.5*-.15) {\bf Figure 4.} A path of $\tilde {\cal P}(i,K)$ in $\tilde K$ (left), 
and the corresponding augmenting\\
   (.53*-.23) path of ${\cal P}(i,K)$ in $K$ (from Figure 3). The shaded region is $\tilde\beta^{-1}\partial\beta{\cal P}_{b(i)}$.\\
    (.86*.81) $x$\\ 
    (.78*.83) $\bar x$\\
\endSetLabels

\AffixLabels{\epsfysize=6.3cm 
\epsfbox{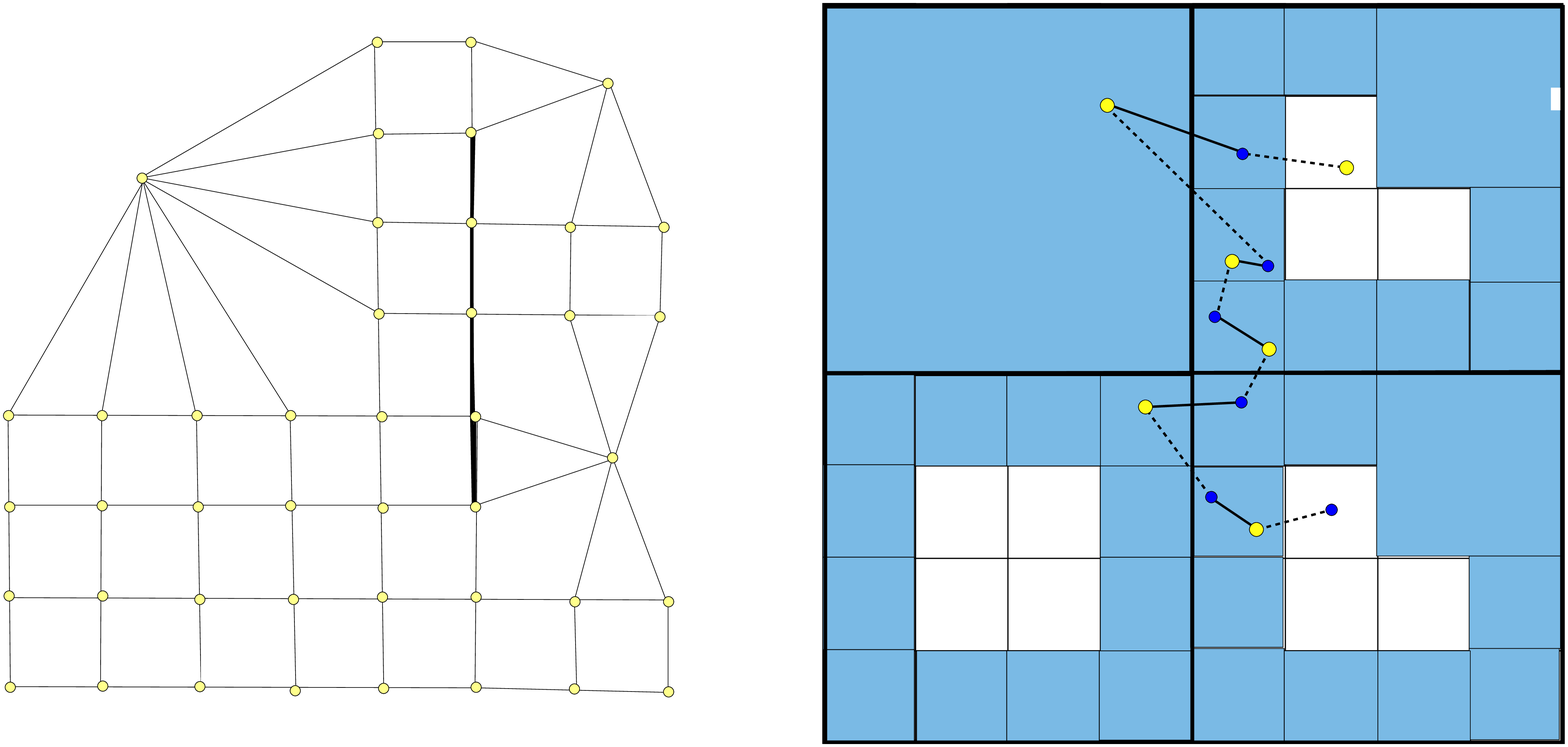}
}

\medskip
\medskip
$$$$

We shall consider the path $P(x,y)$ as an {\it 
augmenting
path}, and replace the edges of $P(x,y)$ in ${\cal M}_{i-1}$ by those
that are not in there. Doing this over all $K\in{\cal P}_{a(i)}$ and $P(x,y)$, we get ${\cal M}_i$. 
For a
particular $K\in {\cal P}_{(a(i))}$, after doing the
``flip" for each augmenting path in ${\cal P}$, all of $K$ but $K\cap 
U_{i}$ is matched. Moreover, the densities of the different edge-lengths
still decay exponentially, see \ref l.gala/.


\def\capa{{\rm cap}}

Let us summarize our conclusion about augmenting paths, also adding a 
claim that the limiting perfect matching exists (which is yet to be 
shown):

\procl p.augmenting
If the conditions (and hence the conclusions) of \ref p.NAGY/ hold, then
{\bf Step 2} is successful, and by the end of stage $i$, all vertices in 
$K\setminus U_i$ are matched. As $i\rightarrow \infty$, we get a
perfect matching in the limit. 
\endprocl

The following lemma tells that the mathings ${\cal M}_j$ stabilize (hence 
proving the last assertion of \ref p.augmenting/),
and 
gives 
the tail probabilities for the limiting matching. For simplicity, in this 
statement and its proof, by an 
edge we mean a pair in one of the ${\cal M}_i$'s. For a vertex $x$, ${\cal M}_i (x)$ will denote the pair of $x$ by ${\cal M}_i$, or the emptyset if $x$ is not matched.
By $||e||$ we denote the 
distance between the two endpoints of $e$.

\procl l.gala
Suppose the conditions in \ref p.NAGY/ are satisfied. 
Then for every vertex $v\in\Z^d$ there is a number $j(v)$ such that $v$ is 
contained 
in an edge $e$ with the property that $e\in {\cal M}_j$, for every $j\geq 
j(v)$. Moreover, we have
$$\P [||e||>r]\leq c \exp (-c' r^{d-2-\epsilon}).$$
\endprocl

\proof
The edges stabilize for the following reason. In order for a vertex $x$ to 
be 
in different edges or no edge in $({\cal M}_j)_j$ infinitely many times, 
it is necessary that either $x$ or ${\cal M}_j (x)$ is contained
in $\cup_{K\in {\cal P}_{a(j)}} {\cal P} (j, K)$, and hence in 
$\tilde\beta^{-1}\partial\tilde\beta\tilde {\cal P}_{b(j)}\cup U_j$, for 
infinitely many $j$'s. The probability of this is 0, since $\sum 
\P [x\in\tilde\beta^{-1}\partial\tilde\beta\tilde {\cal P}_{b(j)}\cup U_j]+\P[{\cal 
M}_j (x)\in\tilde\beta^{-1}\partial\tilde\beta\tilde 
{\cal P}_{b(j)}\cup U_j]\leq 2 \sum
\P [x\in\tilde\beta^{-1}\partial\tilde\beta\tilde {\cal P}_{b(j)}\cup U_j]$ by \ref 
l.parja/, and 
this sum is finite.

Now, let $e$ be as in the claim, and $j(v)$ be the smallest nonnegative 
integer such that $e\in {\cal M}_{j}$ whenever $j\geq j(v)$. If $j(v)=0$, 
then by definition of ${\cal M}_0$, the endpoints of $e$ are in the same 
${\cal P}$-cell, 
and the claim follows by \ref l.M0/. 

Otherwise there is some sequence of edges that $v$ is contained in, until 
it stabilizes from ${\cal M}_{j(v)}$ on. Let $h$ be the greatest number such 
that $v\in U_h$ if such a $U_h$ exists, otherwise $h:=0$. There is an 
$m$, and a sequence $ 
i_1,\ldots,i_m$, $i_1:=h+1$, such that ${\cal 
M}_{i_j}(v)\not ={\cal
M}_{i_j-1}(v)$. We may assume that $m$ is maximal such. Let
the edge containing $v$ in ${\cal
M}_{i_j}$ be called $e_j$. By maximality, $e_m=e$.
Since $v\in U_{i_1 -1}=U_h$, by (II) and (III) we have that 
${\cal M}_{i-1}(v)$ is
in the  
$\tilde\beta^{-1} {\cal P}_{b(i)}$-cell of $v$. 
Thus $||e_1||\leq \sqrt d 2^{b(h)+1}$, unless $v$ is in a ${\cal P}$-cell of 
size $>2^{b(h)}$. This shows the following tail for 
$||e_1||$:
$$\P [||e_1||>r]\leq 2\P [v\in U_j\, ,\, j\geq \log {r\over 2\sqrt d}]+ \P [{\cal 
P}(v) \hbox{ has radius } \geq r/2d^{1/2}]\leq$$
$$\leq c2^{-a([\log r/2\sqrt d])^{1+\epsilon/2}}+
c_0 \exp 
(-c_0'r^{d-2-\epsilon})\leq c \exp (-c'r^{d-2-\epsilon}). \label e.elso
$$
We got the second term on the left of \ref e.elso/ by \ref
l.chernoffnew/, and the first term from 
$$\P[v\in U_i\, , \, j\geq \log r/2\sqrt d]\leq c\sum_{j=[\log r/sqrt d]}^\infty 
\min_{x,y\in U_i} \dist (x,y)^{-d}\leq$$
$$\leq \sum_{j=[\log r/\sqrt d]}^\infty 2^{-a(j)^{1+\epsilon/2}}\leq c 
2^{-a([\log r])^{1+\epsilon/2}},$$
using \ref e.kiemel/.
For the last inequality of \ref e.elso/, we used the choice of $a(i)$.

We get each $e_i$ from $e_{i-1}$ by the ``switch" along some 
augmenting path 
$P(x,y)$, and (by the choice of $h$) we also assumed that they 
are inner 
edges of this path, whenever $i-1>1$. Thus, $e_i$ and 
$e_{i-1}$ being consecutive edges in $P(x,y)$, they share one endpoint, 
while their other endpoints ($x(e_i)$ and $x(e_{i-1}$)
are in adjacent 
${\cal P}$-cells). Hence $\P[||e_i||-||e_{i-1}||\geq r]\leq 
\P[|x(e_j)-x(e_{j-1})|\geq r]\leq 
4\P[{\cal 
P}(o)$ has diameter $\geq r/4]\leq c \exp 
(-c'(r)^{d-2-\epsilon})$, by \ref l.parja/ and \ref l.chernoffnew/. That 
is, 
$$ \P[||e_i||-||e_{i-1}||\geq r]\leq c \exp
(-c'r^{d-2-\epsilon}). \label e.megegy
$$

Also, for each edge $e_i$, if $v$ switched from edge $e_{i-1}$ to edge 
$e_{i}$ in some stage $j$, then in particular, one of the endpoints of 
$e_i$ has to be 
in 
$\tilde\beta^{-1}\partial \tilde\beta {\cal P}_{b(j)}$. We can give the 
following rough bound on this probability: 
$$2\P [o\in\tilde\beta^{-1}(\partial\tilde\beta {\cal P}_{b(j)})]= 2\P 
[\tilde\beta {\cal P}_{b(j)}(o)\in \partial\tilde\beta({\cal 
P}_{b(j)})]\leq$$
$$\leq \P [\diam ({\cal P}(o))\leq b(j)\, ,\, \dist (o, \partial ({\cal 
P}_{b(j)} (o))\leq b(j)]+\P [\diam ({\cal P}(o)\geq b(j)]\leq$$
$$\leq \P [\dist (o, \partial ({\cal
P}_{b(j)} (o)))\leq b(j)]+\P [\diam ({\cal P}(o)\geq b(j)]\leq$$
$$\leq c 2^{-b(j)}b(j)+c\exp (-c' b(j)^{d-2-\epsilon})\leq c\exp (-c' 
b(j)^{d-2-\epsilon}). \label e.bj
$$
Here the bound on the first probability is a consequence of \ref l.MTP/, and the bound on the second probability is by \ref l.chernoffnew/.

Let $A$ be the event that $m>\log R$ and $B$ be the
event that $\bigl\{\{||e_0||>R/\log R\}$ or
$||e_i||-||e_{i-1}||>R/\log R$
for some $i\in \{2,\ldots m\}\bigr\}$. 
We finish the proof by noting that for $||e||>R$ to happen, one has to 
have at least one of $A$ and $B$ hold.
Using \ref e.bj/ to bound the 
probability of $A$ (which is $\leq 2\P 
[o\in \cup_{j=[\log r]}^\infty \tilde\beta^{-1}\partial\tilde \beta ({\cal P}_{b(j)})$),
and \ref 
e.elso/ and \ref e.megegy/ to show $\P[B]\leq c\log R \exp(-c' (R/\log 
R)^{d-2-\epsilon})$,
this gives
$$\P [e>R]\leq c\exp (-c' b(\lceil \log R\rceil )) +c\exp (-c' (R/\log 
R)^{d-2-\epsilon})
\log R .$$
We conclude that
$$\P[e>R]\leq c\exp (-c' R^{d-2-2\epsilon}),$$
since $b(\lceil \log R\rceil )\geq R^d$.
\Qed

So, all what is left is to show the existence of an admissible flow in 
$N(\tilde K)$ of strength $\sigma={|U_{i-1}|+|U_i|\over 2}$ from $B$ to 
$Y$, since then the conditions in \ref p.augmenting/ and 
\ref l.gala/ follow.

By the maxflow-mincut theorem, the existence of such a flow follows if we 
show that every mincut has capacity $\geq \sigma$. 

Before proving this, let us state a rough estimate relating the size of
the boundary of a subgraph of $\tilde K$ (or $K$) to the number of 
dyadic
pseudocubes of a certain type in it. Given $A\subset K$, denote by
$\partial
A$ as before, the set of vertices that have
degree $<2d$ in $A$, and denote by $\partial_j A$ the $j$-neighborhood
of $\partial A$.

\procl l.boundary
Let $\tilde A$ be an induced subgraph of $\tilde K$, and
$A:=\tilde\beta^{-1}
(\tilde A)\subset K$. Then, for the
set $D_\ell (A)$ of dyadic
$k2^\ell$-pseudocubes in $A$ that are not contained in any larger dyadic 
pseudocube in
$A$, we have
$$|D_\ell (A)|\leq 2|\partial A|d^{1/2}(k2^\ell)^{1-d}.$$
\endprocl

\proof
Every element of $D_\ell (A)$ is contained in
$A\cap\partial_{d^{1/2}k2^{\ell +1}} A$, as
can be seen by induction on $i$. The elements of $D_\ell (A)$ are 
disjoint. Thus
dividing
$|A\cap\partial_{d^{1/2}k2^{\ell +1}}A|$ by the
lower bound $(k2^{\ell -1})^d$ on the
volume of an element in $D_\ell$, we get an upper bound $2d^{1/2}|\partial
A|(k2^\ell )^{1-d}$ for $D_\ell (A)$.
Here we were also using $|A\cap\partial_{d^{1/2}k2^{\ell +1}}A|\leq 
|\partial A| d^{1/2} k2^{\ell +1}$.
\Qed

\procl p.flow
There is an admissible flow of strength $\sigma$ from $B$ to $Y$ in 
$N(\tilde K)$.\endprocl

\proof
We will prove by contradiction. So, suppose that there is some  
$\pii$ minimal edge-cutset
between $B$ and 
$Y$, and that $\capa (\pii)<\sigma=\sum_x \capa(\{B,x\})=\sum_x \capa(\{Y,x\}$. Hence $\pii$ contains some 
edge of $\tilde K$. 
Let $\gamma=\pii\cap \tilde K$, and $C_1,C_2,\ldots, C_m$ be the 
connected 
components of $(\tilde K\cap \partial\tilde\beta{\cal 
P}_{b(i)})\setminus 
\gamma$. 
One can find a set 
$\gamma '$ of 
edges of 0 capacity, such that $\tilde K\setminus (\gamma\cup\gamma ')$ 
has components $C_1 ',\ldots, C_m '$, such that that 
$C_i\subset C_i '$, and $\gamma '$ minimal with this property (i.e., for 
any $e\in 
\gamma\cup \gamma '$, $\tilde K\setminus (\gamma\cup\gamma ')\cup \{e\}$
has less than $m$ components). 
Furthermore, one has:
$$|\gamma '|\leq 2^{d(b(i)+1)+1}|\gamma\cap\partial \tilde\beta {\cal 
P}_{b(i)}|. \label e.gamma
$$
This is true because by definition, every cell $C$ of ${\cal P}_{b(i)}$ is 
either a pseudocube of size $2^{b(i)}$ in ${\cal Q}_{b(i)}$, or an 
elementary cell from ${\cal P}$ (which is thus contracted by $\tilde 
\beta$), hence every class of $\tilde\beta {\cal P}_{b(i)}$ has 
cardinality 
$\leq 2^{db(i)+1}$. We can assign to each element of $\gamma\cap C$, 
$C\in\tilde\beta {\cal P}_{b(i)}$, some $e\in \partial\tilde\beta {\cal 
P}_{b(i)}$ in the $\tilde\beta{\cal P}_{b(i)}$-cell of
$C$, to get \ref e.gamma/.

Note that the edges in $\gamma '$ have costs 0 by the definition of 
$N(\tilde K)$.
Hence we may replace $\pi$ by $\pi\cup\gamma '$ and $\gamma$ by $\gamma 
'$, to have (by \ref e.gamma/):

\procl l.gamma
There exists a minimal cutset $\pi$  between $B$ and $Y$ in $N(\tilde K)$, $\gamma=\tilde K\cap\pi$, such 
that the number of components in $\tilde K\setminus\gamma$ is the same as 
the number of components in $(\tilde K\cap \partial\tilde\beta{\cal
P}_{b(i)})\setminus
\gamma$, and further, $|\gamma|\leq (2^{d(b(i)+1)+1}+1) |\gamma\cap \partial
\tilde\beta{\cal
P}_{b(i)}|$.\endprocl

Now, for each $C_j$, one of 
$N_{Y,C_j}$ and $N_{B,C_j}$ (call it $N_j$) has to belong to 
$\pii$, 
otherwise there is a path from $B$ to $Y$ through $C_j$ that avoids 
$\pii$.

\procl l.apro 
There is a $j'\in \Z$ such that
$$\sum_{j\not =j'}|\partial(\tilde\beta^{-1} C_j)|\leq 
c |\tilde\beta^{-1} \gamma|\leq  
c {|\gamma|}\leq c' 2^{b(i)d}|\gamma\cap \partial\tilde\beta {\cal 
P}_{b(i)}|
.$$ \endprocl

\proof
The first inequality is by \ref l.surface/. The second is simply because 
after contracting by $\tilde \beta$ we kept multiple edges.
The last, rough inequality is 
true because of \ref l.gamma/. 
\Qed



%

Apply \ref l.surface/ to $\Gamma:=\tilde\beta^{-1} (\gamma)$ and 
$\Delta_j := \tilde\beta^{-1} (C_j)$. 
We may assume that $i=1$ is the (possible) exception in \ref 
l.surface/ (and equivalently that $C_1$ is the exception in \ref l.apro/), 
and 
assume by symmetry that $N_{1}=N_{B,C_1}$. Furthermore, we may assume that 
every $\Delta_j$ is such that 
$$|\Delta_j\cap U_i|>1, \label e.nagyobb
$$
otherwise if $\Delta_j\cap U_i=\{x\}$, it is easy to see that we could remove part of $\partial 
C_j$ (edges of capacity $k^{{d\over d-1}}/12d'> 1$) from $\pii$, and add 
the edge between $\bar x$ and $B$ or $Y$, to get a cutset of smaller cost 
than $\pii$.
Similar argument works if $\Delta_j\cap U_i=\emptyset$.

As mentioned, if $\pii$ is a minimal cutset, then 
$$\capa (\pii)\leq \capa 
(N_{B,\tilde K}),  \label e.csillag
$$
since $N_{B,\tilde K}$ is a cutset itself.

\def\constant{{C}(b(i))}

Use notation ${C}(b(i)):=(c2^{b(i)d})^{-1}k^{d/(d-1)}/(12d')$, where the first factor and the constant $c$ there is coming from \ref l.apro/, while the second factor is the capacity of edges in $\gamma\cap\partial\tilde\beta {\cal P}_{b(i)}$.
By \ref l.apro/, using that every edge of $\gamma\cap \partial\tilde\beta {\cal P}_{b(i)}$ has capacity $k^{d/(d-1)}/(12d')$, we can bound the left hand 
side of \ref e.csillag/ as  
$$\capa(\pii)=\capa(\gamma)+\sum_{j=1}^m \capa(N_j)\geq 
\constant |\gamma|+\sum_{j=1}^m \capa(N_j)=$$
$$=\constant |\gamma|+\capa(N_1)+\sum_{j=2}^m \capa(N_{B,C_j})-{\bf 
1}_{N_j=N_{Y,C_j}}(\capa(N_{B,C_j})-\capa(N_{Y,C_j}))\geq$$
$$\geq \constant |\gamma|+\capa(N_{B,C_1})+ \sum_{j=2}^m 
\capa(N_{B,C_j})-{\bf
1}_{N_j=N_{Y,C_j}}\sur(\Delta_j\cap U_{i-1})-|\Delta_j\cap 
U_i|.$$

Note that we were using our assumption on $C_1=\tilde\beta (\Delta_1)$ for 
the last inequality.
Putting this fact together with the last inequality and \ref e.csillag/,
we obtain
$$\capa(N_{B,\tilde K})\geq  \constant |\gamma| +\sum_{j=1}^m 
\capa(N_{B,C_j})-\sum_{j=2}^m({\bf
1}_{N_j=N_{Y,C_j}}\sur(\Delta_j\cap U_{i-1})+|\Delta_j\cap 
U_i|).$$
(We define the empty sum as
0. That corresponds to the case when
$\gamma$ is not a cutset for $\tilde K$.)
This implies
$$\sum _{j=2}^m |\Delta_j\cap U_i|+\sum_{j=2}^m \sur (\Delta_j\cap 
U_{i-1})\geq \constant |\gamma|. \label e.sokadik
$$ 
On the other hand, the minimal distance between elements of 
$U_{i-1}$ is $\geq 2^{{a(i-1)(1+\epsilon/2) \over d}}$ by the choice of 
$U_{i-1}$ (see \ref e.kiemel/). Let 
$A(i):= 2^{a(i)(1+\epsilon/2)/d}$. 
We have that
$$\sur (\partial_{A(i-1)}\Delta_j\cap U_{i-1})
\leq |\partial_{A(i-1)}\Delta_j\cap U_{i-1}|
\leq
|\{x\, :\, B_{A(i-i)}(x)\subset\partial_{2A(i-1)}\Delta_j\}\cap 
U_{i-1}\leq$$
$$\leq |\partial_{2A(i-1)}\Delta_j|/A(i-1)^d
\leq c|\partial
\Delta_j|/A(i-1)^{d-1}. \label e.yan
$$
Using notation $D_\ell$ coming from 
\ref l.boundary/, observe that 
$$\Delta_j\subset 
\partial_{A(i-1)} \Delta_j\cup(\cup_{\ell=\log (A(i-1))}^{a(i)}\cup_{X\in 
D_\ell (\Delta_j)} X).$$
Subadditivity $\sur 
(H\cup H')\leq \sur (H)+\sur (H')$ yields
$$\sur (\Delta_j\cap U_{i-1})\leq \sur (\partial_{A(i-1)}\Delta_j\cap 
U_{i-1})+\sum_{\ell=\log (A(i-1))}^{a(i)}\sum_{X\in D_\ell (\Delta_j)} 
\sur 
(X\cap U_{i-1})\leq  \label e.jelol
$$
$$\leq c |\partial \Delta_j| A(i-1)^{1-d}+\sum_{\ell=\log 
(A(i-1))}^{a(i)} |D_\ell (\Delta_j)|f((k2^\ell )^d).$$
Here we used \ref e.yan/ for the first term, and for the second term we 
used the fact that dyadic 
pseudocubes are not bad by definition.

For the first term in \ref e.sokadik/, we get the following bound. In the first 
inequality we use \ref l.combined/ and \ref e.nagyobb/, and in the second 
one we use \ref 
l.surface/ and the choice of $C_1$:
$$\sum_{j=2}^m |\Delta_j\cap U_i|\leq c 2^{-a(i)\epsilon/2}\sum_{j=2}^m 
|\partial\Delta_j|\leq c 2^{-a(i)\epsilon/2} |\gamma|.$$

Plugging this and \ref e.jelol/ 
into 
\ref e.sokadik/, 
\def\ideg{2^{-a(i)\epsilon/2}}
we obtain 
$$C(b(i)) |\gamma| \leq c\ideg |\gamma|+\sum_{j=2}^m
\Bigl(c|\partial\Delta_j|A(i-1)^{1-d}
+\sum_{\ell=\log (A(i-1))}^{a(i)}
|D_\ell (\Delta_j)|f((k2^\ell)^d)\Bigr)\leq$$
$$\leq c\ideg |\gamma|+\sum_{j=2}^m
|\partial\Delta_j|\Bigl(cA(i-1)^{1-d}
+\sum_{\ell=\log (A(i-1))}^{a(i)}
(k2^\ell)^{1-d }f((k2^\ell)^d)\Bigr)\leq$$
$$\leq |\gamma|\Bigl(c+ cA(i-1)^{1-d}
+\sum_{\ell=\log (A(i-1))}^{a(i)}
(k2^\ell)^{1-d} (k2^\ell)^{d-1-\epsilon/2}\Bigr).
\label e.fo
$$
The penultimate inequality was a consequence of \ref l.boundary/ applied 
to each $\Delta_j$, 
while the last line follows from the first two inequalities of \ref 
l.apro/. 

We conclude, by the definition of $f$ and $C(b(i))$, and using $\log 
(A(i-1))\geq a(i-1)/d$:

\procl p.fo 
If $\pii$ is some cutset between $B$ and $Y$, $\gamma=\tilde K\cap\pii$, 
and $\capa (\pii)<\min \{\capa (N_{B,\tilde K}),\capa (N_{Y,\tilde K})\}$, 
then 
$$k^{{d\over d-1}} b(i)^{-d}\leq c+C 
2^{(-1+1/d)a(i-1)(1+\epsilon/2)}+\sum_{\ell=a(i-1)/d}^{a(i)}
(k2^\ell)^{-\epsilon/2},
\label 
e.fontos
$$
where $c$ and $C$ are constants that depend only on $d$.
\endprocl

The second term on the right is bounded by a constant (independently of $k$).
By the definitions of $a(i)$ and $b(i)$, for $i\geq 2$, the inequality 
\ref e.fontos/ fails, so there is
no 
minimal cutset between $B$ and $Y$ different from $N_{B,\tilde K}$ or 
$N_{Y,\tilde K}$, showing that
the desired 
flow on $N(\omega)$ exists. 
If $i=1$, if $k$ was chosen large enough, then
\ref e.fontos/ fails. This finishes the proof of \ref p.flow/.
\Qed

From \ref p.flow/ and \ref p.NAGY/ we have the existence of $ 
{\cal P}(i,K)$, and
conclude by \ref p.augmenting/
that the rematching procedure succeeds.
Combined with \ref l.gala/, this proves \ref t.dim3/.

We mention that once the $a(i)$ and $b(i)$ can be chosen to grow with a 
suitable speed, the crucial inequalities above can be made to be true {\it 
for any} $i\geq 1$. The major difficulty is the start, since $a(0)$ and 
$b(0)$ has to be set 1. This is responsible for the isoperimetric 
considerations and ``tightness" in the above computations.

\procl r.f
Note that $f(x)$ could be chosen any other way so that \ref e.fontos/ 
fails, and the tail for the matching is 
given by \ref 
l.M0/. 
All one has to ensure is that the sum on the right of \ref 
e.fontos/ resulting with this new $f$ is finite when summed up to infinity.
E.g. $f(x):= x^{{d-1 \over d}}(\log x \log\log x)^{-1}$ 
gives the bound $\P [o$ is matched to distance $\geq r]\leq \exp 
({-cr^{d-2} \over\log r\log\log r})$.
We have chosen the slightly weaker bound to make the formulas easier to 
follow.
\endprocl

\medbreak
\noindent {\bf Acknowledgements.}\enspace
I am grateful for Ander Holroyd and G\'abor Pete for helpful discussions.

\startbib
                                                                                
  
\bibitem[A]{A} Aldous, D. (2007) Optimal Flow Through the Disordered 
Lattice {\it Ann. Probab.} {\bf 35}, 397-438.

\bibitem[Ale]{Ale} K. S. Alexander. Percolation and minimal spanning
forests
in infinite graphs. {\it Ann. Probab.}, 23(1):87-104, 1995.

\bibitem[AKT]{AKT} Ajtai, M., Koml\'os, J., \and Tusn\'ady, G. (1984) On 
optimal matchings {\it Combinatorica} {\bf 4}, 259-264.
  

\bibitem[CPPR]{gravi} Chatterjee, S., Peled, P., Peres, Y. \and Romik, D.
(2007) Gravitational allocation to Poisson points, (to appear in {\it
Annals of Math.}).

\bibitem[HL]{HL} Heveling, M. \and Last, G. (2005) Characterization of
Palm measures via bijective point-shifts, {\it Ann. Probab.} {\bf 33},
1698-1715.

\bibitem[HP1]{HPregi} Holroyd, A., \and Peres, Y.
(2003) Trees and matchings from point processes, {\it Electron. Comm. 
Probab.}, {\bf 8}, 17-27.
                                                                          
\bibitem[HP2]{HP} Holroyd, A., \and Peres, Y.
(2005) Extra Heads and Invariant Allocations, 
{\it Ann. Probab.} {\bf 33}, 31-52.

\bibitem[HPPS]{HPPS}
Holroyd, A.E., Pemantle, R., Peres, Y., \and Schramm, O. (2007) Poisson 
Matching 
(preprint).

\bibitem[S]{S}
Soo, T. (2007)
Translation-Invariant Matchings of Coin-Flips on $Z^d$, (preprint).
                   
\bibitem[Ta]{Ta} Talagrand, M. (1994) Matching theorems and empirical 
discrepancy computations using majorizing measures, {\it J. of the AMS}, 
Vol. 7, 
{\bf 2}, 455-537. 

\bibitem[Th]{Th} Thorisson, H. (2000) {\it Coupling, Stationarity, and 
Regeneration}, Probability and its Applications, Springer-Verlag, New 
York. 
                                                             
\bibitem[Ti]{T} 
Tim\'ar, \'A. (2004)  Tree and Grid Factors for General point Processes, 
{\it Elec. Comm. 
in Probab.} {\bf 9}, 53-59.

\bibitem[Y]{Y} Yukich, J.E. (1998) {\it Probability theory of classical 
Euclidean 
optimization}, vol. 1675 of Lecture Notes in Mathematics, Springer-Verlag, 
Berlin.

\endbib 
\bibfile{\jobname}
\def\noop#1{\relax}
\input \jobname.bbl

\filbreak
\begingroup
\eightpoint\sc
\parindent=0pt\baselineskip=10pt

Department of Mathematics,
University of British Columbia,
121-1984 Mathematics Rd.,
Vancouver, BC V6T1Z1, Canada
\emailwww{timar[at]math.ubc.ca}{}
\htmlref{}{http://www.math.ubc.ca/$\sim$timar/}
\endgroup

\bye